\documentclass[11pt]{article}
\newcommand{\documentdate}{11 II 2025}

\usepackage{a4wide,latexsym,amsmath,varioref,graphicx,xcolor,amssymb,diagbox,comment}
\usepackage{tikz}
\usepackage{multirow}

\title{\name: an adaptive pruning-aware gradient method}

\author{Margherita Porcelli\thanks{Dipartimento di Ingegneria Industriale, Universit\`a degli Studi di Firenze, Viale Morgagni 40/44, 50134 Firenze, Italia.  Email: {\tt margherita.porcelli@unifi.it}}
\thanks{ISTI--CNR, Via Moruzzi 1, Pisa, Italia. }
\and Giovanni Seraghiti\thanks{University of Mons, Rue de Houdain 9, 7000 Mons, Belgium. Email: {\tt giovanni.seraghiti@umons.ac.be}}
 \thanks{Dipartimento di Ingegneria Industriale, Universit\`a degli Studi di Firenze, Viale Morgagni 40/44, 50134 Firenze, Italia.}
\and Philippe L. Toint
  \thanks{NAXYS, University of Namur, Namur, Belgium. Email: {\tt philippe.toint@unamur.be}.}
}

\newcommand{\beqn}[1]{\begin{equation}\label{#1}}
\newcommand{\eeqn}{\end{equation}}
\newcommand{\req}[1]{(\ref{#1})}
\newcommand{\ms}{\;\;\;\;}
\newcommand{\tim}[1]{\;\; \mbox{#1} \;\;}

\newtheorem{theorem}{Theorem}[section]
\newtheorem{lemma}[theorem]{Lemma}

\newcommand{\numsection}[1]{\section{#1}\setcounter{equation}{0}}
\newtheorem{corollary}[theorem]{Corollary}

\newcounter{algo}[section]
\renewcommand{\thealgo}{\thesection.\arabic{algo}}
\newcommand{\llem}[2]{\vspace{\baselineskip} 
\noindent\framebox[\textwidth]{\parbox{0.95\textwidth}{
\begin{lemma} \label{#1} \rm #2 \end{lemma} } } \vspace{\baselineskip} }

\newcommand{\algo}[3]{\refstepcounter{algo}
\begin{center}\begin{figure}[htbp]
\framebox[\textwidth]{
\parbox{0.95\textwidth} {\vspace{\topsep}
{\bf Algorithm \thealgo : #2}\label{#1}\\
\vspace*{-\topsep} \mbox{ }\\
{#3} \vspace{\topsep} }}
\end{figure}\end{center}}
\newcommand{\bpr}{{\bf Proof.} \hspace{1.5mm}}
\newcommand{\epr}{\hfill $\Box$ \vspace*{1em}}
\newcommand{\proof}[1]{
\begin{list}{}{
\setlength{\topsep}{0.0pt}
\setlength{\partopsep}{0.0pt}
\setlength{\leftmargin}{0.025\textwidth}
\setlength{\rightmargin}{0.5\leftmargin}
\setlength{\labelwidth}{0.5\leftmargin}
\setlength{\labelsep}{0.25\leftmargin}}
\item \bpr #1 \epr \noindent
\end{list}}
\newcommand{\lthm}[2]{\vspace{\baselineskip} 
\noindent\framebox[\textwidth]{\parbox{0.95\textwidth}{
\begin{theorem} \label{#1} \rm #2 \end{theorem} } } \vspace{\baselineskip} }
\newcommand{\ii}[1]{\{ 1, \ldots, #1 \}}
\newcommand{\iiz}[1]{\{ 0, \ldots, #1 \}}
\newcommand{\iibe}[2]{\{ #1, \ldots, #2 \}}
\newcommand{\calA}{{\cal A}} 
\newcommand{\calD}{{\cal D}}

\newcommand{\calO}{{\cal O}} 
\newcommand{\calR}{{\cal R}}
\newcommand{\calS}{{\cal S}}

\renewcommand{\Re}{\hbox{I\hskip -2pt R}}
\newcommand{\smallRe}{\hbox{\footnotesize I\hskip -2pt R}}
\newcommand{\bigfrac}[2]{\frac{\displaystyle #1}{\displaystyle #2}}
\newcommand{\bigsum}{\displaystyle \sum}

\newcommand{\eqdef}{\stackrel{\rm def}{=}}

\newcommand{\tal}[1]{{\normalsize {\sf #1}}}

\newcommand{\flow}{f_{\rm low}}

\DeclareMathOperator*{\average}{average}
\newcommand{\sign}{\text{sign}}
\newcommand{\argmin}{\text{argmin}}
\newcommand{\name}{\text{\sf prunAdag}}
\newcommand{\noopt}{\text{\sf Relevant only}}

\date{\documentdate}

\begin{document}

\maketitle

\begin{abstract}
A pruning-aware adaptive gradient method is proposed which classifies the variables in two sets before updating them using different strategies. This technique extends the ``relevant/irrelevant" approach of \cite{ding2019global,zimmer2022compression} and allows a posteriori sparsification of the solution of model parameter fitting problems. The new method is proved to be convergent with a global rate of decrease of the averaged gradient's norm  of the form $\calO(\log(k)/\sqrt{k+1})$. Numerical experiments on several applications show that it is competitive.
\end{abstract}

\numsection{Introduction}
This paper deals with first-order Objective-Function-Free Optimization (OFFO) and parameter pruning for optimization problems of the form
\begin{equation}
    \min_{x\in\smallRe^n} f(x)
    \label{problem}
\end{equation}
where $f$ is a smooth function from $\Re^n$ to $\Re$. In particular, our framework can be applied to problems where the objective function $f$ is a loss function depending on the variable $x=(x_1,\dots,x_n)$ defining a model's parameters. For example, in several machine learning applications, the training of the neural network model consists in solving a finite sum minimization problem of the form 
\begin{equation*}
    f(x)=\frac{1}{m}\sum_{i=1}^m \ell_i(x),
    \label{loss}
\end{equation*}
where the components of $x$ correspond to parameters of the network and where both $n$ and $m$ are large, the latter giving the number of samples in the training set and $\ell_i$ being per-sample loss functions. Therefore, we will often refer to the components of the problem variable $x$ as the model/problem parameters. 

OFFO algorithms are methods where the objective function is never computed; instead, they rely only on derivative information, that is on the gradient in the first-order case. A class of OFFO methods, known as adaptive gradient algorithms, gained popularity in the machine learning community, emerging as state-of-the-art techniques to train neural networks. Some examples include Adagrad~\cite{duchi2011adaptive,mcmahan2010adaptive}, Adam~\cite{kingma2014adam}, RMSprop~\cite{tieleman2012lecture}, ADADELTA~\cite{zeiler2012adadelta}. All of these methods share the common characteristic of only relying on current and past gradient information to adaptively determine the step size at each iteration. As Gratton et al. suggested in~\cite{gratton2022first,gratton2023multilevel,gratton2024complexity}, adaptive gradient methods can be interpreted as \textit{trust-region} methods (see \cite{conn2000trust,yuan2015recent} for a comprehensive coverage) in which the radius of the trust-region is computed without evaluating the objective function, which makes them significantly more resistant to noise. Specifically, we propose a new OFFO method based on a modified version of Adagrad in the context of parameter pruning. As the name implies, pruning a model refers to the process of reducing its size and complexity, typically by removing or zeroing certain parameters that are considered unnecessary for its performance. In particular, in the neural network context, the goal of pruning is to improve efficiency, reduce overfitting, and speed up inference or training, without sacrificing much predictive accuracy~\cite{reed1993pruning,lecun1989optimal,zhu2017prune}. We refer to our algorithm as pruning Adagrad (\name).

Pruning emerges as a compression technique for neural networks alternative to matrix or tensor factorization~\cite{zhang2015accelerating,denton2014exploiting,yu2017compressing} or quantization~\cite{courbariaux2016binarized,wang2018two,kim2020position}. Pruning can be performed after training at the cost of re-training the model solely on the remaining parameters~\cite{han2015learning}, but this procedure can sometimes be computationally impractical. Alternatively, one can induce sparsity during training through regularization~\cite{yu2012exploiting,louizos2017learning,alvarez2017compression}; however, this strategy is significantly influenced by the choice of the regularization parameter and the level of sparsity of the solution cannot be altered without re-training the model. To address these limitations, {\em pruning-aware} methods have been developed. They require just one training task and they aim at finding a possibly dense solution, which is robust to pruning, in the sense that the performance of the model is not significantly affected when individual parameters (unstructured pruning)~\cite{han2015learning} or groups of parameters (structured pruning)~\cite{hu2016network} are pruned after training. Consequently, a key concept in pruning-aware methods
consists in the choice of the criteria to determine which parameters or group of parameters can be removed with less impact on the model's performance~\cite{hoefler2021sparsity}. This paper is mainly concerned with unstructured pruning, but our approach can be extended to structured pruning as well.

The majority of the pruning-aware schemes share two common aspects. The first consists in classifying all the parameters, at each iteration, into {\em relevant} and {\em irrelevant} according to specific criteria. Secondly, the method promotes the relevant parameters by update rules which usually involve derivative information, and penalizes the irrelevant ones by diminishing their magnitudes or setting them to zero. This latter strategy may be suboptimal in the context of neural network
training, as it has been shown that the importance of network weights can change dynamically during the training process~\cite{guo2016dynamic,mocanu2018scalable,he2018soft}, meaning that zeroing parameters might decrease the ability to capture these changes. Therefore, a controlled decrease of irrelevant parameters is preferable. At the end of the training phase, a model trained by a pruning-aware method has relevant components with larger magnitudes than the irrelevant ones. Finally, irrelevant components are pruned by removing/zeroing those parameters that are below some threshold such that the model matches any desired level of sparsity. We now give an overview of existing pruning-aware schemes and then present our contributions in this context.

\noindent
\paragraph*{Related works.} We distinguish two classes of pruning-aware methods. The first divides parameters into relevant and irrelevant and then updates them following various rules. This approach is referred to as {\em activation selection} in~\cite{ding2019global}. The second updates parameters all at once and forces the magnitude of irrelevant parameters to decrease by adding specific sparsity-inducing constraints to the problem, see~\cite{pokutta2020deep,zimmer2022compression}. 

The first pruning-aware approach is presented in~\cite{ding2019global,molchanov2016pruning,theis2018faster} and employs Taylor series to measure the
importance of a parameter by estimating the impact of its removal/zeroing on the value of the objective function in \req{problem}. More specifically, let $x_k=(x_{1,k},\dots,x_{n,k}) \in \mathbb{R}^n$ be the $k$-th iterate of the method, then in the first-order case, the relevance of the parameter $x_{i,k}$ at iteration $k$ is measured by
$$f(x_{1,k},\dots, \underbrace{ 0}_{x_{i,k}},\dots,x_{n,k})-f(x_{1,k},\dots,x_{i,k},\dots,x_{n,k})=\frac{\partial f}{\partial x_i}(x_k)(0-x_{i,k})+o(x_{i,k}^2).$$ 
Then, the $T_k$ parameters with largest values of $|\frac{\partial f}{\partial x_i}(x_k)(0-x_{i,k})|$ are classified as relevant, as they are the ones that most significantly affect the objective function's value. Using this criterion, Ding et al.~\cite{ding2019global} propose to optimize relevant components via momentum SGD, while gradually decreasing the magnitudes of all others classified as irrelevant. One of the advantages of this approach is that the number of relevant components $T_k$ can be chosen at each iteration, in order to match a prescribed level of sparsity. An adaptive choice to select $T_k$, using the $\ell_0$-norm of the parameters, is proposed in~\cite{ding2019global}. Although very efficient in practice, the convergence of gradient descent methods using Taylor series approach to classify relevant and irrelevant parameters is not analyzed in~\cite{ding2019global}. Furthermore, it is not clear how much the irrelevant components can be reduced at each iteration without affecting the convergence of the method towards a stationary point.

A second class of pruning-aware methods consists in adding sparsity-inducing constraints to the formulation in \req{problem} and solving the deriving constrained optimization problem using the (stochastic) Frank-Wolfe (SFW) algorithm~\cite{zimmer2022compression,lu2022learning}. Specifically, the constraints considered are $T$-sparse polytope and $T$-support-norm-ball for unstructured pruning and group-$T$-support-norm-ball for structured pruning, see~\cite{argyriou2012sparse,rao2017group,pokutta2020deep}. Since we are primarily interested in unstructured pruning in a deterministic setting, we briefly describe the FW two-step framework~\cite{frank1956algorithm,levitin1966constrained} for solving a $T$-support-norm-ball constrained problem, which is used for comparison in Section~\ref{sec:numer}. Let $T>0$ and $\tau \in \mathbb{R}$, then the $T$-support-norm-ball $\mathcal{C}_{T}(\tau)$ is defined as
$$\mathcal{C}_{T}(\tau)=\text{conv}\{x\in \mathbb{R}^n \ |\ \lVert x \rVert_0 \leq T, \quad  \lVert x \rVert_2 \leq \tau\},$$
where conv$(\cdot)$ denotes the convex hull. The FW algorithm is applied to the sparse-constrained model 
$$\min_{x\in \mathcal{C}_{T}(\tau)} f(x)$$
and consists of two main steps. First, a descent direction is computed, solving the linear minimization oracle (LMO)
$$
v_k = \ \argmin_{v \in \mathcal{C}_T(\tau)} \langle v, g_k \rangle,
$$
where $g_k$ is the gradient of $f$ at the $k$-th iterate $x_k$.  
The optimal solution of this problem is given by
\begin{equation}\label{vk}
v_{i,k}  =
  \left\{ \begin{array}{ll}
 -\tau \  g_{i,k} /\| g_{k} \|_{\calR_k} & \text{if }
  i \in \mathcal{R}_k,\\
  0 & \text{otherwise,}
  \end{array}\right.
  \end{equation}
  for $i=1, \dots,n$, where $\|\cdot\|_{\calR_k}$ is the Euclidean norm of the subvector containing the indices of the $T$ largest components of the absolute value of $g_k$. Next, the new feasible iterate is obtained as a convex combination of the past iterate and the descent direction, that is
  \begin{equation}\label{learn}
  \ x_{k+1}=x_k+\eta_k( v_k-x_k),  
  \end{equation}
for some learning rate $\eta_k \in [0,1]$. In machine learning or neural networks applications, SFW is often employed and utilizes the stochastic gradient instead of the actual gradient $g_k$. Specifically, the convergence of SFW for a finite sum minimization problem in the nonconvex setting has been established in~\cite{reddi2016stochastic} and it has been extended to the gradient rescaling version by Zimmer et al. in~\cite{zimmer2022compression}.
The SFW algorithm has demonstrated performance comparable to state-of-the-art pruning methods \cite{zimmer2022compression}; however, like the deterministic version, it has certain drawbacks. Specifically, the stepsize $\tau$ is related to the radius of the constrained region, which must be predetermined and cannot vary throughout the iterations; as well as the number of components $T$ that are updated using gradient information. Moreover, the method is extremely sensitive to the tuning of its parameters such as the scalar $\tau$, the learning rate $\eta_k$, and the number of gradient components $T$ in the search direction. 

\medskip
Our \name\ method is a new pruning-aware scheme. At each iteration, we classify as relevant those parameters corresponding to the $T$ largest directional derivative in magnitude, similarly to Zimmer et al. in~\cite{zimmer2022compression}, but the overall classification and updating strategy are significantly different.
\noindent
\begin{enumerate} 
    \item  We develop a new adaptive strategy to separately update parameters that extends the distinction between relevant and irrelevant ones. Specifically, we introduce the concepts of {\em optimisable} and {\em decreasable} parameters. We consider all parameters that benefit from being updated using derivative information as optimisable, including all relevant parameters and those irrelevant parameters that can be both optimized and penalized simultaneously, using gradient information. We define as decreasable those irrelevant parameters that are not penalized by a gradient update and therefore require a specific penalization strategy to decrease their magnitudes.    
    \item We propose to use the Adagrad step for updating the optimisable parameters, while we develop an Adagrad-like trust-region framework to gradually decrease the magnitude of the decreasable ones without relying on gradient information or any function evaluations. 
    \item We prove the convergence of \name\ in the deterministic case with complexity of $O(\log(k) / \sqrt{ k+1})$.
    \item We validate our method on several preliminary examples obtained from signal processing, dictionary learning and binary classification applications.
\end{enumerate}\medskip

\paragraph*{Organization of the paper}
Our presentation is structured as follows. Section \ref{method} introduces our \name\ method in Algorithm \ref{alg:PA} and discusses both the criterion used for classifying parameters at each iteration into optimisable and decreasable and the different update rules for the two classes of parameters. Convergence of the method is analyzed in Subsection  \ref{subsec:conv}. Section \ref{sec:numer} 
presents a set of illustrative examples of the behaviour of \name\ when applied to
random least-squares, sparse signal recovery, sparse coding step in dictionary learning, and binary classification with logistic loss function. A brief conclusion is provided in Section~\ref{sec:concl}.

\noindent
\paragraph*{Notations.} Throughout the paper we adopt the following notations. At the $k$-th iteration, the gradient $g$ of the function $f$ evaluated at the current iterate $x_k$ is denoted as $g_k = g(x_k)$ and its $i$-th component by $g_{i,k}$. Moreover, the superscript $T$ denotes the
transpose and $v_{i,k}$ denotes the $i$-th component of a vector
$v_k\in \Re^n$. Unless specified otherwise, $\|\cdot\|$ is the
Euclidean norm on $\Re^n$ and $\|x\|_{\mathcal{I}}=\| x_{i \in \mathcal{I}} \|$ is the Euclidean norm when we consider only the indices in $\mathcal{I}$. We denote as $A^C$ the complementary set of $A$ in $\ii{n}$.  We use the notation $\lceil x \rceil$ for the minimum integer greater than $x$. Given two sequences $\{\alpha_k\}$ 
and $\{\beta_k\}$ of non-negative reals, we also say that $\alpha_k$ is $\calO(\beta_k)$ is there exists
a finite constant $\kappa$ such that $\lim_{k\rightarrow\infty} (\alpha_k/\beta_k) \leq \kappa$.

\numsection{A first-order sparsity inducing adaptive gradient method without regularization}
\label{method}

In the following discussion, we make the standard assumptions for first-order complexity analysis, which ensure that problem \req{problem} is consistent.
\vspace*{2mm}
\noindent
\begin{description}
\item[AS.1:] the objective function $f(x)$ is continuously differentiable;
\item[AS.2:] its gradient $g(x)$ is Lipschitz continuous with
   Lipschitz constant $L\geq 0$, that is
   \[
   \|g(x)-g(y)\| \le L \|x-y\|
   \]
   for all $x,y\in \Re^n$; 
\item[AS.3:] there exists a constant $\flow$ such that, for all $x$, $f(x)\ge \flow$.
\end{description}

We now motivate and describe the \name\ method. It falls in the class of pruning-aware methods that, at a given iteration, first classifies the variables/parameters in relevant and irrelevant ones before updating them.

\noindent
Following Zimmer et al. in~\cite{zimmer2022compression,pokutta2020deep} and given a target $T$ for the cardinality of the relevant parameters, we first define the set $\calR_k$ of $\calR$elevant parameters as the set of indices corresponding to the $T$ largest directional derivatives in magnitude:
\beqn{calR-def}
\calR_k \eqdef \{ i \in \ii{n} \mid |g_{i,k}| \mbox{ is one of the $T$ largest components of }|g| \}.
\eeqn
The parameters indexed by $\calR_k^C$ are thus considered irrelevant.  The idea is then to optimize on the relevant parameters and to reduce the magnitude of the others.  In our approach, optimization is performed by applying the component-wise version of Adagrad, in which the step $s_{i,k}$ in the $i$-th variable is given by
\beqn{adag-step}
s_{i,k} = -\frac{g_{i,k}}{w_{i,k}^{\calO}},
\eeqn
for some suitably chosen weight $w_{i,k}^{\calO}$ derived from the history of past gradients.
However, it may happen that the step \req{adag-step} that would be taken by the optimizer for some of the irrelevant parameters does actually reduce their magnitudes: this happens when the signs of $x_{i,k}$ and $g_{i,k}$ coincide. Thus it makes sense to take the optimization step on these parameters as well, provided it remains within reasonable bounds compared to $|x_{i,k}|$. 
More formally, we deem the $i$-th component ($i\in \calR_k^C)$ to be ``optimisable" (in the sense
that its update can be performed by \req{adag-step}) if it belongs to the set $\calA_k$ of $\calA$cceptable indices
\beqn{calA-def}
\calA_k 
\eqdef \left\{ i\in \calR_k^C \mid {\rm sign}(g_{i,k}) = {\rm sign}(x_{i,k}) \tim{and}
              a_{i,k} \leq \left|\frac{g_{i,k}}{w_{i,k}^{\calO}}\right|\leq b_{i,k}  \right\},
\eeqn
where $a_{i,k}$ and $b_{i,k}$ are suitable bounds. For example, one can choose $a_{i,k}$ as a fraction of the absolute value of $x_{i,k}$ (further practical examples are given below). As a consequence, the index set of the optimizable parameters at iteration $k$ is given by 
\beqn{calO-def}
\calO_k \eqdef \calR_k \cup \calA_k.
\eeqn
We also see that the ``decreasable" parameters, whose magnitude we wish to $\calD$ecrease by other means than \req{adag-step}, have indices in the set $\calD_k=\calR_k^C\setminus\calA_k$.

While the course of action for the ``optimisable" parameters is clear (apply \req{adag-step}), what to do with the other parameters is, at this stage, less obvious. Our proposal is to extend
our Adagrad-based approach by defining a step for the parameters in $\calD_k$ which is also of the form (optimality measure / weight), where now the optimality measure is no longer given by gradients values, but by the magnitude of the irrelevant parameters themselves (as we wish to drive them to zero, if possible), that is a step of the form
\beqn{irr-step}
s_{i,k}^{\calD} = -\frac{x_{i,k}}{w_{i,k}^{\calD}},
\eeqn
where the weight $w_{i,k}^{\calD}$ is now derived from the history of past irrelevant values of the $i$-th parameter.

We are now ready to state the \name\ algorithm more formally in Algorithm \ref{alg:PA}.

\algo{alg:PA}{\tal{\name}}
{
\begin{description}
\item[Step 0: Initialization. ]
  A starting point $x_0$, a target number of relevant parameters $T\in \iibe{2}{n}$, a constant $\varsigma\in(0,1)$, and two initial weight vectors $w_{i,-1}^\calO=\varsigma$ and $w_{i,-1}^\calD=\varsigma,$ for $i=1,\dots,n$ are given. Set $k=0$. 

\item[Step 1: Select relevant parameters.] Compute the gradient $g_k $ and define
 \beqn{Ik_top}
    \calR_k= \{i \in \{1,\dots,n\} \ | \ |g_{i,k}| \ \mbox{is one of $T$ largest components of} \ |g_k| \}
  \eeqn
  \item[Step 3: Optimisation weights.] Compute 
  \beqn{temp_weig}
  w_{i,k}^\calO = \sqrt{(w_{i,k-1}^\calO)^2+g_{i,k}^2}
      \quad (i \in \ii{n}) 
   \eeqn
   \item[Step 4: Classify components.] Define $a_{i,k}, b_{i,k} \geq 0$ for $i \in \calR_k^C$,  and
   \beqn{calA-def2}
   \calA_k=\left\{i \in \calR_k^C \mid \sign(x_{i,k})=\sign(g_{i,k}) \tim{and} a_{i,k} \leq \left\lvert \frac{g_{i,k}}{w_{i,k}^{\calO}}\right\rvert \leq b_{i,k}  \right\}.
   \eeqn
   Then set $\calO_k = \calR_k\cup\calA_k$ and $\calD_k = \calO_k^C$.
  \item[Step 5: Compute a step for the optimisable components.] Compute
  \beqn{adagrad weights}
  s_{i,k}=-\frac{g_{i,k}}{w_{i,k}^\calO} \quad (i \in \calO_k).
  \eeqn

\item[Step 6: Compute a step for the decreasable components.] 

   Compute 
 \beqn{bottom_weights}
        w_{i,k}^\calO = w_{i,k-1}^\calO, \ms
        w_{i,k}^\calD = \sqrt{(w_{i,k-1}^\calD)^2+x_{i,k}^2}, \quad 
        s_{i,k}^L=-\frac{x_{i,k}}{w_{i,k}^\calD} \quad (i\in\calD_k) 
    \eeqn
    
    \beqn{bottom_weights_fixed}
        w_{i,k}^\calD = w_{i,k-1}^\calD \quad (i \in \calO_k).
    \eeqn
Compute a step $s_{i,k}$ for all the decreasable components in $\calD_k$ such that 
\beqn{dis_abs_S}
|s_{i,k}| \leq |s_{i,k}^L| \qquad (i\in \calD_k)
\eeqn
\beqn{dec_cond}
\sum_{i \in \calD_k} g_{i,k} s_{i,k} \leq 0.
\eeqn
 
\item[Step 7: New iterate.] Define
   $$  \qquad x_{k+1} = x_k + s_k,$$ 
    increment $k$ by one and return to Step~1.
\end{description}
}

A few comments are useful after this formal description.
\begin{enumerate}
\item The weights defined in \req{temp_weig} and \req{bottom_weights} superficially look identical to weights used by Adagrad. There is however a crucial difference for our purpose: each of these updating formula only selects in the past those iterations for which the considered component (the $i$-th) is either optimisable (for the former) or decreasable (for the latter).  More specifically, if we set
\beqn{gOxD-def}
g_{i,k}^\calO=\left\{\begin{array}{cc}
    g_{i,k} & i \in \calO_k \\
    0 & \text{otherwise,}
\end{array}\right. \qquad x_{i,k}^\calD=\left\{\begin{array}{cc}
    x_{i,k} & i \in \calD_k \\
    0 & \text{otherwise,}
\end{array}\right.
\eeqn
then one verifies, using the first part of \req{bottom_weights} and \req{bottom_weights_fixed}, that
\beqn{weights2}
w_{i,k}^\calO = \sqrt{\varsigma + \sum_{j=0}^k (g_{i,k}^\calO)^2}
\tim{ and }
w_{i,k}^\calD = \sqrt{\varsigma + \sum_{j=0}^k (x_{i,k}^\calD)^2}.
\eeqn
\item We have left the choice of the exact technique to define $s_k^\calD$ very general, as long as \req{dis_abs_S}-\req{dec_cond} hold (these conditions can be interpreted as ``trust-region" conditions with radius $|s_{i,k}^L|$). One simple technique is, for example, to set
\begin{equation}
    s_{i,k}^\calD=\left\{ \begin{array}{ll}
        -\sign(x_{i,k})  \cdot\min[a_{i,k},|s_{i,k}^L|]  & i \in \calS_k \eqdef  \{i \in \calD_k \ | \ \sign(x_{i,k})=\sign(g_{i,k})\}   \\
         0 & \mbox{otherwise.}
    \end{array}\right.
\end{equation}

\item The specification of the set $\calA_k$ is not strictly necessary in our convergence theory below, but our experience indicates that extending the set of optimisable parameters from $\calR_k$ to $\calR_k\cup\calA_k$ is beneficial in practice, as illustrated in Section~\ref{sec:OvsR}.
In what follows, we will consider the choices defined in Table~\ref{tab:versions}. 

\begin{table}[ht]
    \centering
    \begin{tabular}{|l|c|c|}
    \hline
         &  $a_{i,k}$ & $b_{i,k}$\\
    \hline     
     \rule[-3ex]{0pt}{7ex} Version 1  & $\left\lvert\bigfrac{x_{i,k}}{k+1} \right\rvert \cdot \bigfrac{\lVert g_k \rVert_{\calR_k}}{\lVert x_k \rVert_{\calS_k}}$ & $+\infty$ \\
     \hline
     \rule[-2.5ex]{0pt}{7ex} Version 2  & $\left\lvert\bigfrac{x_{i,k}}{k+1}\right\rvert$ & $+\infty$ \\
      \hline
     \rule[-3ex]{0pt}{7ex} Version 3  & $\left\lvert\bigfrac{x_{i,k}}{k+1} \right\rvert \cdot \bigfrac{\lVert g_k \rVert_{\calR_k}}{\lVert x_k \rVert_{\calS_k}}$ & $\lvert x_{i,k} \rvert$ \\
      \hline
     \rule[-2.5ex]{0pt}{7ex} Version 4  & $\left\lvert\bigfrac{x_{i,k}}{k+1}\right\rvert$ & $\lvert x_{i,k} \rvert$ \\
     \hline
    \end{tabular}
    \caption{Four versions of \name, depending on the choice of the bounding sequences $a_{i,k}$ and $b_{i,k}$ used in Step~4 of Algorithm \ref{alg:PA}.}
    \label{tab:versions}
\end{table}

In all four versions, the lower bounding sequence $a_{i,k}$ represents a fraction of the absolute value of $i$-th component of the iterate. In Version~2 and Version~4, the sequence linearly decreases with the iteration number, while in Version~1 and Version~3, in addition to the linearly decreasing factor, we rescale the sequence to match the magnitude of the gradient $g_k$ in $\calR_k$. A similar rescaling of the learning rate has proven effective in SFW~\cite{pokutta2020deep,zimmer2022compression}. In our experience, choosing a lower bound sequence $a_{i,k}$ that decreases as $1/k$ helps to make the convergence of the method faster and reduces the oscillations due to the non-monotone behaviour of the algorithm. On the other hand, the upper bounding sequence $b_{i,k}$ in Version~3 and Version~4 is set to be equal to the iterate magnitude $|x_{i,k}|$. This choice prevents the Adagrad-like step from exceeding $|x_{i,k}|$, thus avoiding sign changes for the components in $\calA_k$. Allowing potential sign changes in some components can be obtained by setting $b_{i,k}=\infty$ (as in Version~1 and Version~2), allowing for a larger set $\calA_k$, as Figure~\ref{Fig:act_evolution} shows. This implies that the number of optimisable components is then larger and the algorithm's speed is potentially enhanced. Conversely, the number of indices in $\calD_k$ is smaller, potentially leading to less effective pruning.  We illustrate the evaluation of the sets $\calO_k$, $\calA_k$ and $\calD_k$ for the four versions of \name\ in Section~\ref{sec:act_set}.
\end{enumerate}

\subsection{Convergence analysis} \label{subsec:conv}

First, and although the boundness of the gradient is not required in our analysis, we assume that the sequence produced by \name\ is bounded, that is, 
 \begin{description}
  \item[AS.4:] The sequence produced by Algorithm \ref{alg:PA} is bounded, i.e. for some $\kappa_x > 0$, $| x_{i,k} | \leq \kappa_x$ for every $i = \{1,\dots,n\}$ and all $k\geq 0$.
\end{description}
The main purpose of this assumption is to ensure that the magnitudes of the decreasable variables/parameters remain finite and thus that the task of reducing them is realistic.
We also assume, without loss of generality, that $\varsigma \leq \left(\frac{8nL}{3}\right)^2$.

\medskip

\noindent
We start the convergence analysis of the \name\ algorithm by stating 
a lemma characterizing the descent properties of the method. 

\llem{lemma:dec_ver3}{
Suppose that AS.1 and AS.2  hold. Then we have that, for all
$j\ge0$,
\beqn{gen-decr}
f(x_{j+1})
\le f(x_j) -\sum_{i=1}^n \frac{ (g_{i,j}^\calO)^2}{ w_{i,j}^\calO}
     + \frac{L}{2} \sum_{i=1}^n \frac{ (g_{i,j}^\calO)^2}{(w_{i,j}^\calO)^2} 
     + \frac{L}{2} \sum_{i=1}^n \frac{ (x_{i,j}^\calD)^2}{(w_{i,j}^\calD)^2}\\
\eeqn
and
\beqn{gen-dec-k}
f(x_0)-f(x_{k+1}) \geq 
\sum_{j=0}^k \sum_{i=1}^n \frac{(g_{i,j}^\calO)^2}{ w_{i,j}^\calO}
   - \frac{L}{2} \sum_{j=0}^k \sum_{i=1}^n \frac{ (g_{i,j}^\calO)^2}{(w_{i,j}^\calO)^2} 
   - \frac{L}{2} \sum_{j=0}^k \sum_{i=1}^n \frac{ (x_{i,j}^\calD)^2}{(w_{i,j}^\calD)^2}.
\eeqn
}

\proof{
Using \req{dec_cond}, we have that, at each iteration $j$ of Algorithm \ref{alg:PA}, 
\beqn{cond_s_bot}
 g_j^T s_j = \sum_{i\in\calO_j} g_{i,j}s_{i,j} + \sum_{i\in\calD_j} g_{i,j}s_{i,j}
 \leq \sum_{i\in\calO_j} g_{i,j}s_{i,j} <0.
 %(g_j^\calO)^T s_j^\calO + (g_j^\calD)^T s_j^\calD 
 %\leq (g_j^\calO)^T s_j^\calO<0.
  \eeqn
  Therefore, using Assumptions AS.1 and AS.2, \req{cond_s_bot}, \req{adagrad weights}, \req{bottom_weights}, and \req{dec_cond} we derive that
\beqn{f_dec}
\begin{aligned}
f(x_{j+1})&\leq f(x_j)+g_j^T s_j +  \frac{L}{2} ||s_j||^2 \\
%&\leq f(x_j)+ (g_j^\calO)^T s_j^\calO +  \frac{L}{2} ||s_j^\calO||^2 +  \frac{L}{2} ||s_j^\calD||^2 \\
&\leq f(x_j)+ \sum_{i\in\calO_j} g_{i,j}s_{i,j} +  \frac{L}{2} \sum_{i\in \calO_j}s_{i,j}^2 +  \frac{L}{2} \sum_{i\in \calD_j}s_{i,j}^2 \\
&\leq f(x_j) -\sum_{i \in \calO_j} \frac{ g_{i,j}^2}{ w_{i,j}^\calO}
      + \frac{L}{2}  \sum_{i \in \calO_j} \frac{ g_{i,j}^2}{ (w_{i,j}^\calO)^2} + \frac{L}{2} \sum_{i \in \calD_j}  \frac{ x_{i,j}^2}{(w_{i,j}^\calD)^2}\\
&=f(x_j) -\sum_{i = 1}^n \frac{ (g_{i,j}^\calO)^2}{ w_{i,j}^\calO}
     + \frac{L}{2} \sum_{i=1}^n \frac{ (g_{i,j}^\calO)^2}{ (w_{i,j}^\calO)^2} 
     + \frac{L}{2} \sum_{i=1}^n \frac{ (x_{i,j}^\calD)^2}{(w_{i,j}^\calD)^2}      
\end{aligned}
\eeqn
and \req{gen-decr} therefore holds. Finally, summing for $j=0,\dots,k$ gives (\ref{gen-dec-k}).
}

\noindent
We next proceed by recalling the following lemma from \cite{WuWardBott18}, which is crucial in the derivation of upper bounds for the second and third terms in equation \req{gen-dec-k}.

\llem{gen:series}{Let $\{a_k\}_{k\ge 0}$ be a non-negative sequence, $\xi>0$ and define, for each $k \geq 0$,
$b_k = \sum_{j=0}^k a_j$.  Then

\beqn{alphasup1series-bound}
\sum_{j=0}^k  \frac{a_j}{(\xi+b_j)}
\le  \log\left(\frac{\xi + b_k}{\xi} \right).
\eeqn
}

\noindent
We now state our main complexity result, partly inspired by \cite[Theorem~3.2]{gratton2024complexity}.

\lthm{theorem:allmu}{Suppose that AS.1--AS.4 hold. Assume that $\calR_k$ defined in \req{Ik_top}, has cardinality $T$ for every $k$ and that the
\name\ algorithm is applied to problem \req{problem}.
If we define
$
\Gamma_0 \eqdef f(x_0)-\flow,
$
then,
\beqn{gradbound}
\average_{j\in\iiz{k}}\|g_j\|^2 \le \lceil n/T \rceil \frac{\theta(k)}{k+1}
\eeqn
with \beqn{k6-def}
\theta(k) \eqdef \max\left\{
\varsigma,\bigfrac{\varsigma}{2} e^ \frac{\Gamma_0}{ n L},
32n^2L^2
\,\left|W_{-1}\left(-\frac{\sqrt{\varsigma}}{8nL}\right)\right|^2, \ 2 \left( \Gamma_0 +nL\log\left( 1 + \frac{(k+1)\kappa_x^2}{\varsigma} \right) \right)^2
\right\},
\eeqn
where $W_{-1}$ is the second branch of the Lambert function \cite{Corletal96}.
}

\proof{
Let us first observe that bounding the average of  $\|g_j^\calO\|^2$ allows us to derive a bound for the average of the norm of the actual gradient $\|g_j\|^2$. Indeed, since $\calR_j\subseteq \calO_j$ contains the largest components of the gradient at iteration $j$ and its cardinality is always equal to $T$, we have that
\begin{equation*}
    \|g_j\|^2 
    \leq  \lceil n/T \rceil \sum_{i \in \calR_j} g_{i,j}^2 
    \leq \lceil n/T \rceil \sum_{i \in \calO_j} g_{i,j}^2 
    =  \lceil n/T \rceil \|g_j^\calO\|^2
\end{equation*}
and therefore
\begin{equation}
  \average_{j\in\iiz{k}}\|g_j\|^2 \leq \lceil n/T \rceil\average_{j\in\iiz{k}} \|g_j^\calO\|^2.  
  \label{mod_grad_avg}
\end{equation}

Now, using (\ref{gen-dec-k}) and the fact that the sequence $w_{i,k}^\calO$ in \req{adagrad weights} is increasing in $k$ for every $i$, we have that
\begin{equation*}
    f(x_{k+1}) \leq f(x_0)- \sum_{j=0}^k \frac{ \|g_{j}^\calO\|^2}{ \max_{i\in\{1,\dots,n\}} w_{i,k}^\calO}
    + \frac{L}{2} \sum_{j=0}^k \sum_{i=1}^n \frac{ (g_{i,j}^\calO)^2}{(w_{i,j}^\calO)^2} 
    + \frac{L}{2} \sum_{j=0}^k \sum_{i=1}^n \frac{ (x_{i,j}^\calD)^2}{(w_{i,j}^\calD)^2},
\end{equation*}
from which we obtain that
\beqn{dec_grad}
\sum_{j=0}^k \frac{ \|g_{j}^\calO\|^2}{ \max_{i\in\{1,\dots,n\}} w_{i,k}^\calO} 
\leq \Gamma_0 +\frac{L}{2} \sum_{j=0}^k \sum_{i=1}^n \frac{ (x_{i,j}^\calD)^2}{(w_{i,j}^\calD)^2}
              +\frac{L}{2} \sum_{j=0}^k \sum_{i=1}^n \frac{ (g_{i,j}^\calO)^2}{(w_{i,j}^\calO)^2}.
\eeqn
We then use Lemma \ref{gen:series} with $a_{j}=(g_{i,j}^\calO)^2$ and the first part of \req{weights2}, yielding
\beqn{use_lemma}
\begin{aligned}
\sum_{i =1}^n \sum_{j=0}^k  \frac{ (g_{i,j}^\calO)^2}{(w_{i,j}^\calO)^2}
&= \sum_{i =1}^n \sum_{j=0}^k  \frac{ (g_{i,j}^\calO)^2}{\varsigma+\sum_{\ell=0}^j (g_{i,\ell}^\calO)^2} 
\leq \sum_{i =1}^n\log \left( \frac{1}{\varsigma} \left( \varsigma + \sum_{\ell=0}^k (g_{i,\ell}^\calO)^2  \right) \right)\\
& \leq n \log\left( 1 + \frac{1}{\varsigma}  \sum_{\ell=0}^k \|g_{\ell}^\calO\|^2 \right).
\end{aligned}
\eeqn
Using again Lemma \ref{gen:series}, this time with $a_{j}=(x_{i,j}^\calD)^2$, the second part of \req{weights2} and Assumption AS.4, we deduce that
\beqn{use_lemma_x}
\begin{aligned}
\sum_{i =1}^n \sum_{j=0}^k  \frac{ (x_{i,j}^\calD)^2}{ (w_{i,j}^\calD)^2}
&= \sum_{i =1}^n \sum_{j=0}^k  \frac{ (x_{i,j}^\calD)^2}{\varsigma+\sum_{\ell=0}^j (x_{i,\ell}^\calD)^2} 
\leq \sum_{i =1}^n\log \left( \frac{1}{\varsigma} \left( \varsigma + \sum_{\ell=0}^k (x_{i,\ell}^\calD)^2  \right) \right)\\
& \leq n \log\left( 1 + \frac{(k+1)\kappa_x^2}{\varsigma}   \right).
\end{aligned}
\eeqn
Combining (\ref{dec_grad}), (\ref{use_lemma}), (\ref{use_lemma_x}) therefore gives that
\beqn{final_grad}
\sum_{j=0}^k \frac{ \|g_{j}^\calO\|^2}{ \max_{i\in\{1,\dots,n\}} w_{i,k}^\calO} \leq \Gamma_0+\frac{n L}{2}  \log\left( 1 + \frac{(k+1)\kappa_x^2}{\varsigma} \right)+\frac{n L}{2} \log\left( 1 + \frac{1}{\varsigma}  \sum_{j=0}^k \|g_{j}^\calO\|^2 \right),
\eeqn
where the first logarithmic term depends on the upper bound on the entries of the iterate in Assumption AS.4, while the second depends on the sum of the norms of the past optimisable gradients.

We now proceed by analyzing two separate cases.

\textbf{Case 1.} Assume first that the contribution given by the optimisable gradients exceeds that of the first logarithmic term in $k$, that is 
\begin{equation*}
   \label{case_1_log}
     \log\left( 1 + \frac{1}{\varsigma}  \sum_{j=0}^k \|g_{j}^\calO\|^2 \right) \geq \log\left( 1 + \frac{(k+1)\kappa_x^2}{\varsigma} \right). 
\end{equation*}
Then, the inequality in \req{final_grad} becomes 
\beqn{final_grad_case1}
\sum_{j=0}^k \frac{ \|g_{j}^\calO\|^2}{ \max_{i\in\{1,\dots,n\}} w_{i,k}^\calO} \leq \Gamma_0 +n L \log\left( 1 + \frac{1}{\varsigma}  \sum_{j=0}^k \|g_{j}^\calO\|^2 \right).
\eeqn
Assume now that
\beqn{case_1}
 \sum_{j=0}^k \|g_{j}^\calO\|^2 
 \geq \max \left[ \varsigma, \bigfrac{\varsigma}{2} e^\frac{\Gamma_0}{ n L} \right],
\eeqn
which implies
\begin{equation*}\label{eqtmp}
1+\frac{1}{\varsigma} \sum_{j=0}^k \|g_{j}^\calO\|^2 \leq \frac{2}{\varsigma} \sum_{j=0}^k \|g_{j}^\calO\|^2 \quad \text{and} \quad  \Gamma_0 \leq n L \log\left(  \frac{2}{\varsigma}  \sum_{j=0}^k \|g_{j}^\calO\|^2 \right),
\end{equation*}
and observe from the first part of \req{weights2} that, for all
$i\in \{1, \dots, n\}$ and all $k$,
\beqn{bound_weights}
w_{i,k}^\calO \leq \sqrt{\varsigma+\sum_{j=0}^k \| g_j^\calO\|^2}.
\eeqn
Thus, from \req{final_grad_case1}, \req{case_1}, and (\ref{bound_weights})  we obtain that
$$%\beqn{}
\frac{ \sum_{j=0}^k\|g_{j}^\calO\|^2}{\sqrt{2\sum_{j=0}^k \|g_{j}^\calO\|^2 }} \leq 2 n L \log\left(  \frac{2}{\varsigma}  \sum_{j=0}^k \|g_{j}^\calO\|^2 \right),
$$%\eeqn
that is
\beqn{equation_log}
\frac{\sqrt{\varsigma}}{2} \sqrt{\frac{2}{\varsigma} \sum_{j=0}^k \|g_{j}^\calO\|^2}  \leq 4 n L \log\left(  \sqrt{\frac{2}{\varsigma}  \sum_{j=0}^k \|g_{j}^\calO\|^2 }\right).
\eeqn
If we now define
\begin{equation}\label{param_eq}
\gamma_1=\frac{\sqrt{\varsigma}}{2}, \qquad \gamma_2=4nL, \qquad u= \sqrt{\frac{2}{\varsigma}  \sum_{j=0}^k \|g_{j}^\calO\|^2 },
\end{equation}
we first note that our assumption that $\varsigma \leq (\frac{8nL}{3})^2$ ensures that $\gamma_2 > 3\gamma_1$. Furthermore, the inequality \req{equation_log} can then be rewritten as
\beqn{tosolvemuhalf}
\gamma_1 u \le \gamma_2 \log(u).
\eeqn
Let us denote by $\psi(u) \eqdef \gamma_1 u - \gamma_2 \log(u)$. Since
$\gamma_2 > 3 \gamma_1$, the equation $\psi(u)=0$ admits two roots $u_1 \leq
u_2$ and \req{tosolvemuhalf} holds for $u\in[u_1,u_2]$.
The definition of $u_2$ then gives that
\[
\log(u_2)- \frac{\gamma_1}{\gamma_2}u_2 = 0,
\]
which is
\[
u_2e^{-\frac{\gamma_1}{\gamma_2}u_2} = 1.
\]
Setting $z = -\frac{\gamma_1}{\gamma_2}u_2$, we obtain that
\[
z e^z = -\frac{\gamma_1}{\gamma_2}.
\]
Thus $z = W_{-1}(-\frac{\gamma_1}{\gamma_2})<0$, where $W_{-1}$ is the second
branch of the Lambert function defined over $[-\frac{1}{e}, 0)$.
As $-\frac{\gamma_1}{\gamma_2} \geq -\frac{1}{3} $, $z$ is well defined and thus
\[
u_2
= -\frac{\gamma_2}{\gamma_1}\,z
= -\frac{\gamma_2}{\gamma_1}\,W_{-1}\left(-\frac{\gamma_1}{\gamma_2}\right)>0.
\]
Therefore, using \req{param_eq}, 
\[
\bigsum_{j=0}^k\|g_{j}^\calO\|^2
= \frac{\varsigma}{2}\,u_2^2
=32n^2L^2
\,\left|W_{-1}\left(-\frac{\sqrt{\varsigma}}{8nL}\right)\right|^2.
\]
Hence, taking the average gives that
\begin{equation}
\average_{j\in\iiz{k}} \|g_{j}^\calO\|^2 \leq 32n^2L^2
\,\left|W_{-1}\left(-\frac{\sqrt{\varsigma}}{8nL}\right)\right|^2 \cdot \frac{1}{k+1}.
\label{avg_case1_1}
\end{equation}
Suppose now that \req{case_1} fails. Then, obviously,
\begin{equation}
\average_{j\in\iiz{k}}  \|g_{j}^\calO\|^2 \leq \max \left[ \varsigma, \bigfrac{\varsigma}{2} e^ \frac{\Gamma_0}{ n L} \right] \cdot \frac{1}{k+1}.
\label{avg_case1_2}
\end{equation} 

\textbf{Case 2.} Consider now the case where 
\begin{equation*}
    \label{case_2_log}
     \log\left( 1 + \frac{1}{\varsigma}  \sum_{j=0}^k \|g_{j}^\calO\|^2 \right) \leq \log\left( 1 + \frac{(k+1)\kappa_x^2}{\varsigma} \right).
\end{equation*}
Then inequality in \req{final_grad} becomes
\beqn{final_grad_case2}
\sum_{j=0}^k \frac{ \|g_{j}^\calO\|^2}{ \max_{i\in\{1,\dots,n\}} w_{i,k}^\calO} \leq \Gamma_0 +nL\log\left( 1 + \frac{(k+1)\kappa_x^2}{\varsigma} \right).
\eeqn
If, on one hand,
\beqn{case_1_2}
 \sum_{j=0}^k \|g_{j}^\calO\|^2 \geq  \varsigma,
\eeqn
then we have that
\begin{equation*}
\bigfrac{1}{\sqrt{2}}\frac{ \sum_{j=0}^k\|g_{j}^\calO\|^2}{\sqrt{\sum_{j=0}^k \|g_{j}^\calO\|^2 }} \leq \Gamma_0 +nL\log\left( 1 + \frac{(k+1)\kappa_x^2}{\varsigma} \right),
\end{equation*}
that is
\beqn{}
\sum_{j=0}^k\|g_{j}^\calO\|^2 \leq 2 \left( \Gamma_0 +nL\log\left( 1 + \frac{(k+1)\kappa_x^2}{\varsigma} \right) \right)^2.
\eeqn
Taking the average then gives that
\beqn{average_1}
\average_{j\in \{0,\dots,k\}} \|g_{j}^\calO\|^2 \leq 2 \left( \Gamma_0 +nL\log\left( 1 + \frac{(k+1)\kappa^2}{\varsigma} \right) \right)^2 \cdot \frac{1}{k+1}.
\eeqn
If, on the other hand, \req{case_1_2} does not hold, then
\beqn{average_2}
\average_{j\in \{0,\dots,k\}} \|g_{j}^\calO\|^2 \leq \frac{\varsigma}{k+1}. 
\eeqn
We finally deduce \req{gradbound} by considering the largest upper among \req{avg_case1_1}, \req{avg_case1_2}, \req{average_1}, \req{average_2} and using  \req{mod_grad_avg}. 
}

\noindent
Theorem \ref{theorem:allmu} demonstrates that \name\ has complexity of $O(\log(k)/\sqrt{k+1})$, in contrast with the original Adagrad algorithm for which the average gradient norm decreases like $O(1/\sqrt{k+1})$. Therefore, \name\ can in general be (marginally) slower than the original Adagrad algorithm. This is due to the fact that deviating the gradient flow to reduce the magnitude of decreasable parameters comes at a complexity cost. Indeed, Lemma \ref{lemma:dec_ver3} shows that the decrease of the function at each iteration, expressed by equation \req{gen-dec-k} is governed by two positive sums in the right-hand side, depending on the gradient of optimisable components $g_{i,k}^\calO$ and on the magnitudes of decreasable components $x_{i,k}^\calD$. By contrast, the descent lemma (Lemma 2.1 in~\cite{gratton2022first}) of first-order OFFO methods, such as Adagrad, only contains a single positive sum depending on the gradient entries. It is easy to see from equation \req{gen-dec-k} that the function value eventually decreases monotonically if the two positive sums in the right-hand side are dominated by the first, negative sum. However, we might expect a monotonic decrease of the function for large enough weights for the Adagrad algorithm, while the non-monotonic behaviour may persist until convergence for \name. This is because the magnitude of the decreasable components in $\calD_k$ is (hopefully) small at convergence but often non zero for all of them. Therefore, the second positive sum in \req{gen-decr}, depending on the decreasable components in $\calD_k$, might remain significant,
even for large $k$, but its growth is fortunately bounded by a term in $\calO(\log(k))$. 
\begin{itemize}
    \item If we start from an initial point that is far from any stationary point of the problem, we might expect the largest components of the gradient to be large in magnitude for several iterations, potentially exceeding the magnitudes of the decreasable components of the iterate. This is exactly the scenario described in Case~1 of the proof of Theorem \ref{theorem:allmu}. Consequently, we observe an empirical Adagrad-like almost linear decrease, as suggested by equations \req{avg_case1_1} and \req{avg_case1_2}, until the contribution of the decreasable components $x_{i,k}^\calD$ exceeds that of the optimisable gradient components $g_{i,k}^\calO$. 
    \item Theorem \ref{theorem:allmu} proves that the average norm of the gradient converges to zero. However, we cannot expect the same behaviour for the decreasable components of the iterates, which are only likely to be small in magnitude at convergence. This means that after a certain iteration, the contribution of the optimisable gradient will be smaller than that of the decreasable components in equation \req{dec_grad}. Thus, Case~2 in the proof typically occurs for large $k$. As a consequence, we may then expect, for a sufficiently large $k$, a decrease of the order of $\log(k)/\sqrt{k+1}$ as suggested by equation \req{average_1}.
\end{itemize}
\noindent
Both these observations suggest that we might expect a faster decrease during the first iterations, followed by a potential slowdown when $\|g_k^\calO\|$ becomes small, as can be observed for the non-rescaled Versions 2 and 4 in Figure~\ref{Fig:rand_grad_spar} and Figure~\ref{Fig:sparco_grad} in the next section. In general, the rescaling of the step in the decreasable components promotes a faster convergence and this behaviour is not observed. Nevertheless, a fast convergence is not the only purpose of a pruning-aware method since robustness to pruning also needs to be considered.

\numsection{Numerical experiments}\label{sec:numer}

We now present numerical tests on a variety of problems from different applications originating in
\begin{itemize}
\item a standard class of randomly generated under-determined linear least-squares problems,
\item the SPARCO library for sparse signal recovery~\cite{van2007sparco} (as supplied by S2MPJ \cite{GratToin24}), which contains test cases from signal processing applications specifically designed for sparse optimization,
\item the ``sparse coding step" in dictionary learning problem,
\item minimizing the logistic function in binary classification problems on several well-know data sets. 
\end{itemize}
These test problems were chosen to test if \name\ is able to enhance convergence to a solution which is robust to pruning. (It is not our purpose to compare \name\ to state-of-the-art techniques in each of the applications considered.)

We implemented four versions of \name\ 
with different choices of  $a_{i,k}$ and $b_{i,k}$ at Step 4 of Algorithm \ref{alg:PA}. We considered the rules Version 1, Version 2, Version 3 and Version 4 described in Table~\ref{tab:versions}, yielding the corresponding implementations \name-V1, \name-V2, \name-V3\ and \name-V4. We also implemented the deterministic version of the Frank-Wolfe method of \cite{zimmer2022compression} for unstructured pruning, considering two learning rates $\eta_k$ in \req{learn}: a linearly decreasing rate $\eta_k=1/(k+1)$ (FW1) and the adaptively rescaled rate~\cite{pokutta2020deep} 
$\eta_k=\min\left[\frac{\beta \lVert g_k\rVert_{ \calR_k}}{\lVert v_k-x_k \rVert},1\right]$, 
where $\calR_k$ is the set containing the indices associated to the $T$ largest components of the gradient in absolute value, $v_k$ it is given in \req{vk} and $\beta \in (0,1)$ (FW2).
Finally, \name\ reduces to the standard Adagrad algorithm if one chooses $T=n$ and avoids performing any classification of the parameters. We set $\varsigma=1/100$ as for \name\ and the Adagrad implementations.

The algorithms are implemented in MATLAB  R2021b on a  64-bit Samsung/Galaxy with 11th Gen Intel(R) Core(TM) i5-1135G7 @ 2.40GHz  and 8 GB of RAM, under Windows 11 version 23H2.

All experiments are randomly initialized with a normalized starting point with exactly $T$ nonzero entries, representing a feasible point for FW1 and FW2 whenever the $T$-support-norm-ball has radius greater than 1. We performed 20 runs per test problem and report averaged results in the Appendix, while we show representative plots on a single run in what follows. We set $T= n/10$ in both our method and the FW variants (but we propose an analysis on the impact of the parameter $T$ in subsection~\ref{subsec:inf_T}).  
We also tuned, for each set of experiments, the stepsize $\tau$ in the FW implementations (named $\tau_1$ and $\tau_2$ in FW1 and FW2, respectively) and the learning parameter $\beta$ of FW2 by trial and error.  Each algorithm is terminated when either
$\|g_k\| \le 10^{-9}$ or when $10^4$ iterations have been performed, unless otherwise specified. 

Throughout this section we perform pruning on the solution \textit{after optimisation} using two different strategies, depending on the analysis considered. Therefore, either we fix a scalar threshold $\delta$ and we remove all the components below this threshold, or we choose a sparsity level of the solution $\sigma$ (in percentage) and we set the threshold $\delta$ to achieve a $\sigma$-sparse solution after pruning. (Note that the sparsity $\sigma$ of the solution after pruning  and the percentage of parameters removed are equivalent quantities; therefore, we will use the two terms interchangeably.)

We now describe two performance measures that will be used below to discuss the numerical results. Given a fixed $\delta>0$, let $x$ be the approximated solution returned by some pruning-aware algorithm and denote by $\Bar{x}$ the pruned version of $x$, that is the vector $x$ whose components with magnitude less than $\delta$ have been zeroed out.
We then define a quality measure to evaluate robustness of the pruning approach considering
\begin{equation}
    \rho \eqdef \lVert g(\Bar{x}) \rVert
    \tim{ and }
    \omega \eqdef  \sqrt{| f(\Bar{x})- f(x) |}.
    \label{qual_meas}
\end{equation}
Both these quantities provide estimates of how much pruning perturbs the solution from $x$. A small value for $\rho$ means that the pruned solution $\Bar{x}$ is close to stationarity, while a small $\omega$ means that the objective function's value at the pruned solution does not differ much from that at $x$. 

\subsection{Random linear least-squares}
In the first set of experiments, we consider a class of randomly generated linear least-squares of the form
\begin{equation}
        f(x)= \frac{1}{2} \|Ax-b\|^2,  
         \label{ls}
    \end{equation} 
with five different matrices $A \in \Re^{m \times n}$ as in~\cite{wen2010fast} and~\cite{porcelli2014variable}. Specifically, we choose $A$ such that
    \begin{itemize}
        \item[$A$1)] $A$ is randomly generated from a normal Gaussian distribution,
        \item[$A$2)] $A$ is a random orthogonal matrix,
        \item[$A$3)] $A$ is random and has orthogonal columns,
        \item[$A$4)] $A$ is random and has orthogonal rows,
        \item[$A5$)] $A \in (0,1)$ generated from a Bernoulli distribution,
        \item[$A6$)] $A$ is obtained from the discrete cosine transform matrix of dimension $n$.
    \end{itemize}
We set $m=100$ and $n=1000$, making the problem under-determined with a relatively large-dimensional subspace of solutions. We then generate a random solution $x^*=\text{randn}(n,1)$ and compute the right-hand side as $b=Ax^*$. 

\subsubsection{Optimisable vs. relevant}\label{sec:OvsR}

\begin{figure}[htb]
\centering
\begin{minipage}[b]{.49\linewidth}
  \centering
  \centerline{\includegraphics[width=\linewidth]{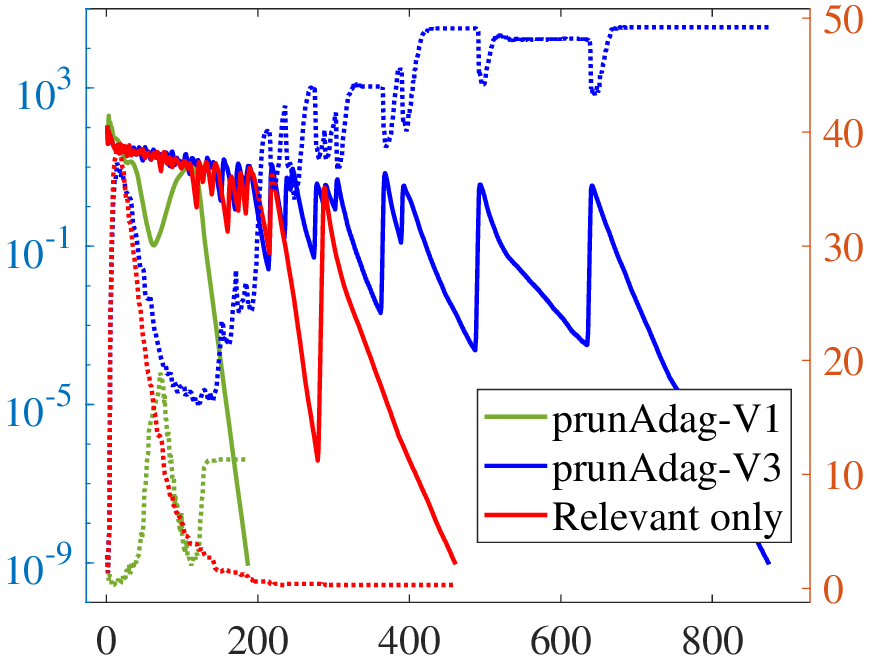}}
  \centerline{(a) gradient norm (left), sparsity (right)}\medskip
\end{minipage}
\hfill
\begin{minipage}[b]{0.49\linewidth}
  \centering
  \centerline{\includegraphics[width=\linewidth]{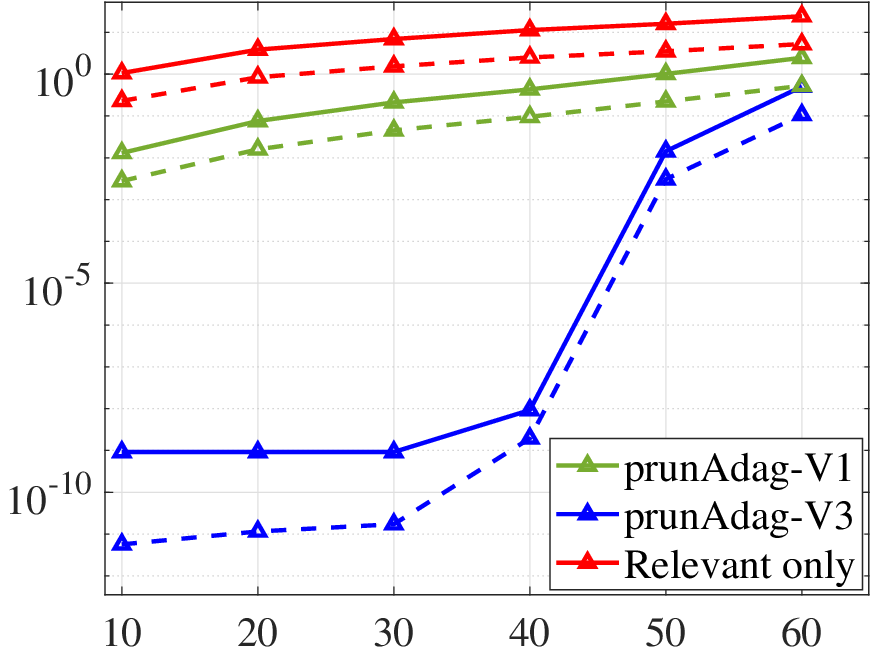}}
  \centerline{(b) measures $\rho$ (solid) and $\omega$ (dashed)} \medskip
\end{minipage}

\centering
\begin{minipage}[b]{.49\linewidth}
  \centering
  \centerline{\includegraphics[width=\linewidth]{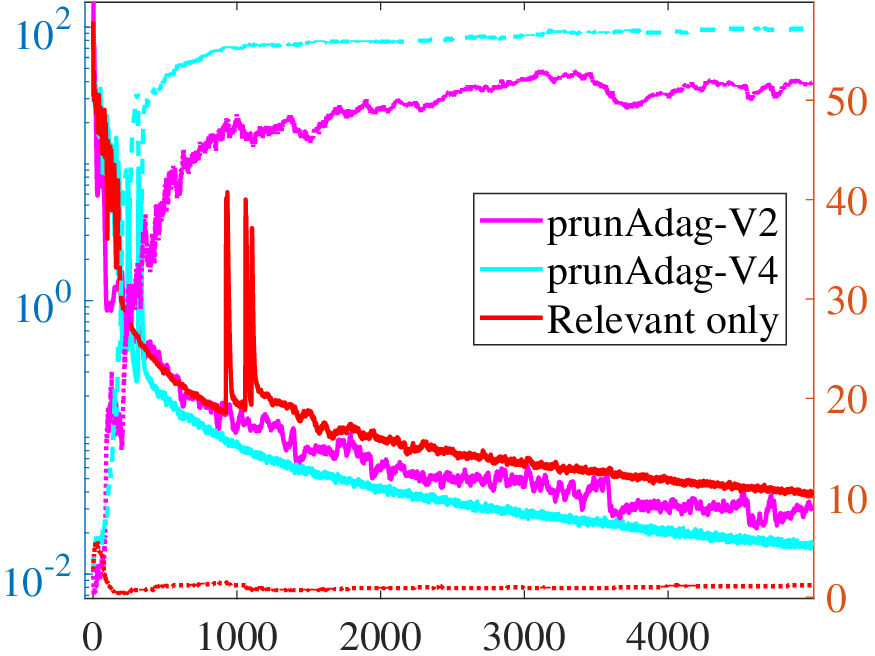}}
  \centerline{(a) gradient norm (left), sparsity (right)}\medskip
\end{minipage}
\hfill
\begin{minipage}[b]{0.49\linewidth}
  \centering
  \centerline{\includegraphics[width=\linewidth]{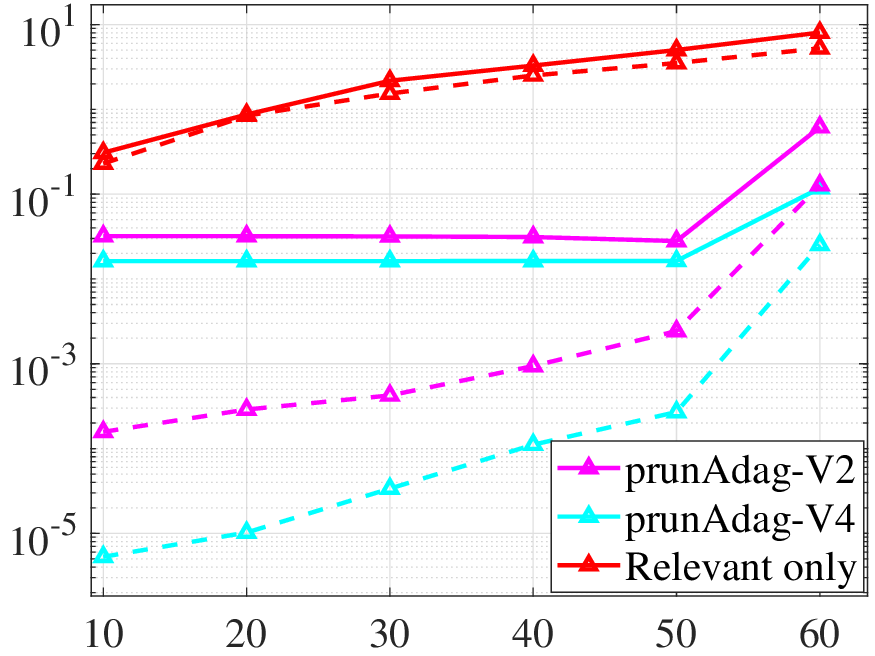}}
  \centerline{(b) measures $\rho$ (solid) and $\omega$ (dashed)} \medskip
\end{minipage}
\caption{Effect of introducing the class of optimisable parameters within the minimization framework of \name. On the right, (a) gradient norm (solid line) on the left $y$-axis and percentage of parameters below $\delta=10^{-3}$ (dotted line) on the right $y$-axes along the iterations. (b) Error measures $\rho$ (continuous) and $\omega$ (dashed) for different percentages of pruned components after the optimization.(Random least-squares A3)}\label{Fig:ex_no_opt_v1_v3}
\end{figure}

We first use the random least-squares problem A3 to show the advantage of extending the parameter classification from relevant/irrelevant to optimisable/decreasable. Figure~\ref{Fig:ex_no_opt_v1_v3} shows the results of running \name-V1, \name-V2, \name-V3 and \name-V4 (which differ by the choice of the bounding sequences $a_{i,k}$ and $b_{i,k}$ as defined in Table~\ref{tab:versions}) in comparison with a version, denoted \noopt, which is identical to \name, except that Step~4 is reduced to the definitions $\calO_k = \calR_k$ and $\calD_k = \calR_k^C$. Thus no "acceptable" parameter is added to the list of relevant ones in \noopt.

\noindent 
The contrast between \noopt\ and the other version is striking in terms of achieved number of components below the sparsity threshold and, consequently, of robustness to pruning. Not only \noopt\ is much less efficient in this respect, but it does not converge faster to a stationary point than the other versions (except when compared to \name-V3, which achieves the best sparsity). In our experience, this behaviour is quite general and, in our view, fully justifies the introduction of $\calA_k$ in Step~4.

\clearpage

\subsubsection{Classification and convergence for the four \name\ versions}\label{sec:act_set}

We next illustrate the evolution of the index sets $\calO_k$, $\calA_k$ and $\calD_k$ along iterations for the four versions of \name. 

\begin{figure}[ht]
\centering
\begin{minipage}[b]{.49\linewidth}
  \centering
  \centerline{\includegraphics[width=\linewidth]{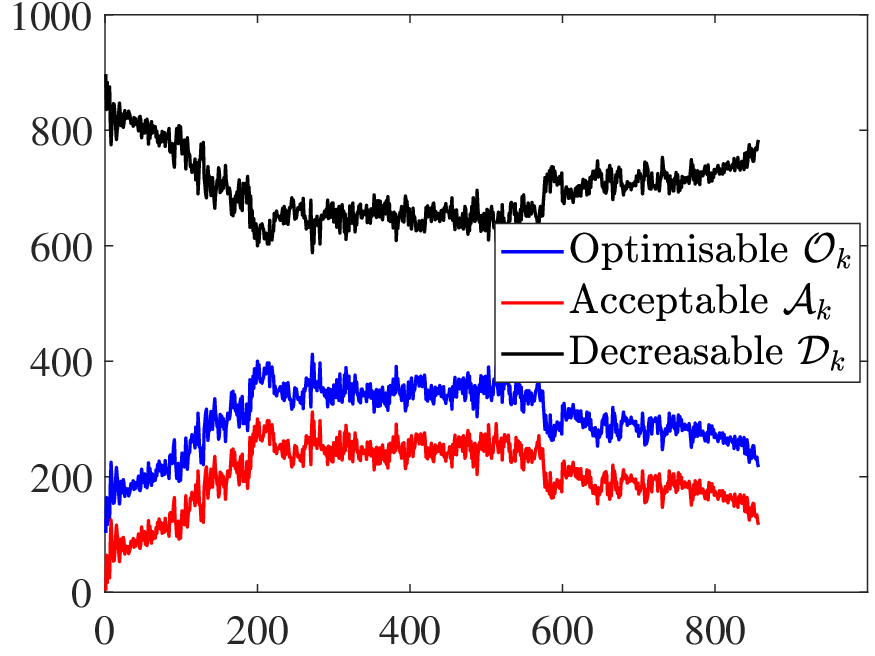}}
  \centerline{(a) \name-V1}\medskip
\end{minipage}
\hfill
\begin{minipage}[b]{0.49\linewidth}
  \centering
  \centerline{\includegraphics[width=\linewidth]{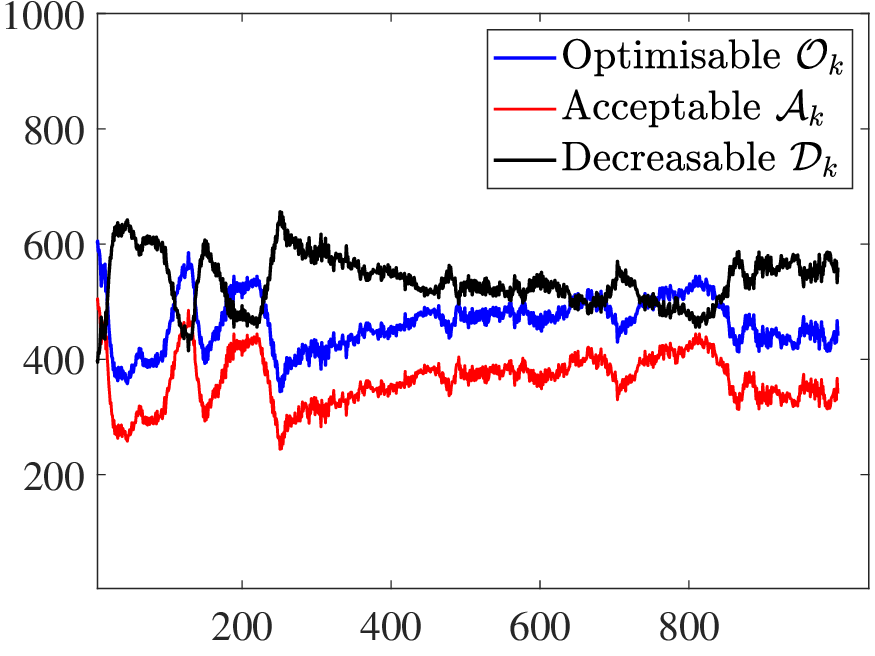}}
  \centerline{(b) \name-V2} \medskip
\end{minipage}

\begin{minipage}[b]{0.49\linewidth}
  \centering
  \centerline{\includegraphics[width=\linewidth]{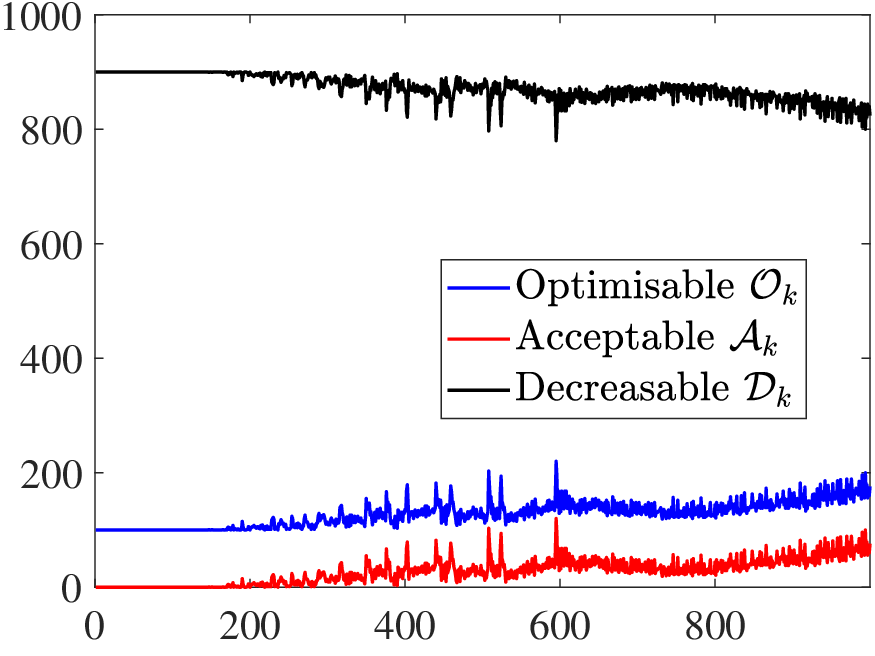}}
  \centerline{(c) \name-V3} \medskip
\end{minipage}
\begin{minipage}[b]{0.49\linewidth}
  \centering
  \centerline{\includegraphics[width=\linewidth]{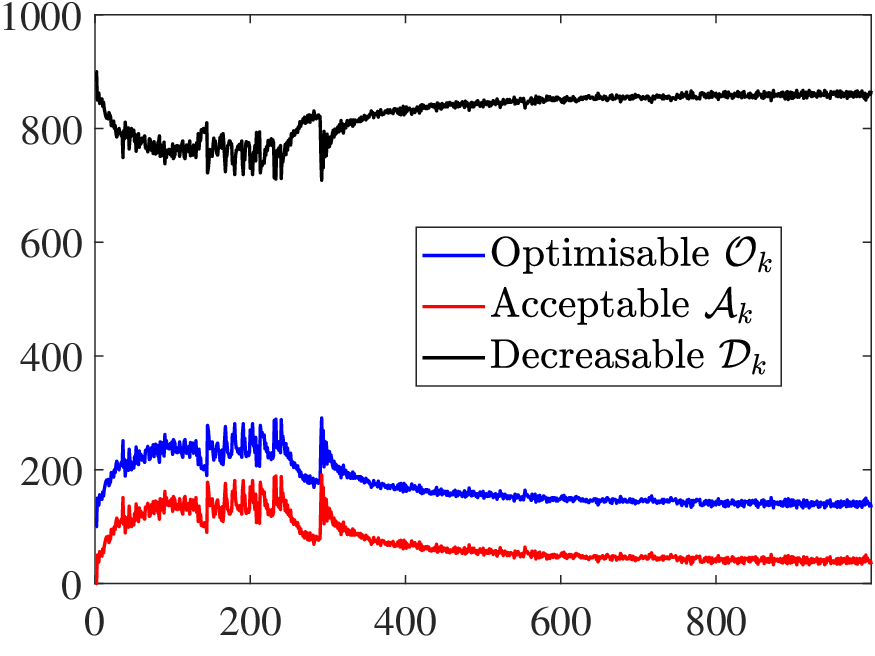}}
  \centerline{(d) \name-V4} \medskip
\end{minipage}
\caption{Dynamic of parameters' classification in the sets $\mathcal{O}_k$, $\mathcal{A}_k$, and $\mathcal{D}_k$ at each iteration. (Random least-squares A1)}
\label{Fig:act_evolution}
\end{figure}

Figure \ref{Fig:act_evolution} shows the dynamic of the indices classified in $\calO_k$, $\calA_k$, 
and $\calD_k$ for all the versions of \name\, when solving the random least-squares problem A1. As shown by this example, this dynamic does vary significantly from version to version,
and we may expect these variations to affect performance.  We observe that robustness to pruning is unsurprisingly better when the size of $\calD_k$ is larger, favoring \name-V3 and \name-V4 on average. If we now turn to speed of convergence to stationary points, the
conclusion is less clear-cut. While, for \name-V1 and \name-V3 (that are the methods using gradient rescaling in $a_{i,k}$), the speed seems to improve with the size of $\calA_k$, the effect is more problem-dependent for the (unscaled) \name-V2\ and \name-V4.

\subsubsection{Influence of $T$}
\label{subsec:inf_T}
\begin{figure}[ht]
\centering
\begin{minipage}[b]{.49\linewidth}
  \centering
  \centerline{\includegraphics[width=\linewidth]{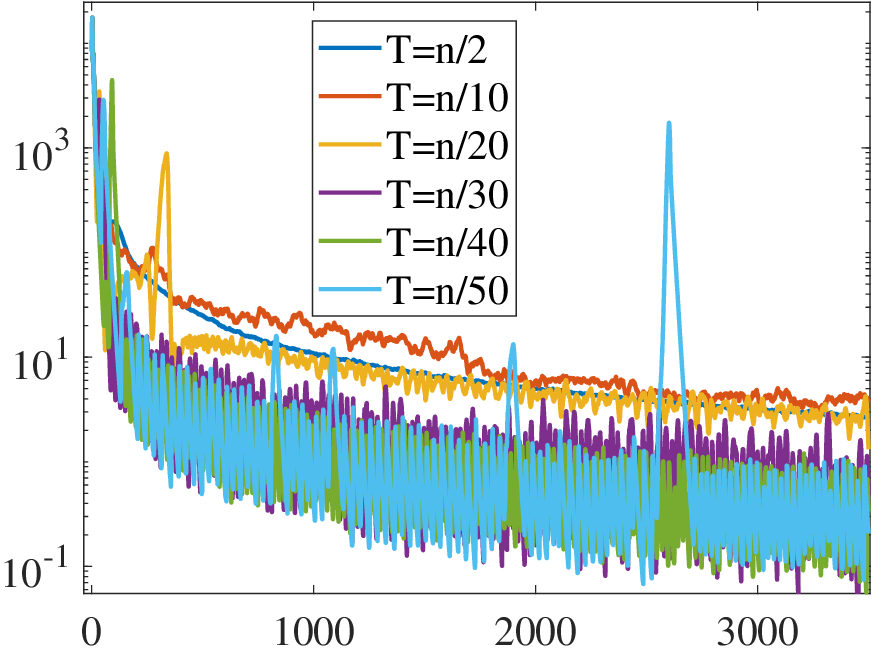}}
  \centerline{(a) \name-V2 - gradient norm}\medskip
\end{minipage}
\hfill
\begin{minipage}[b]{0.49\linewidth}
  \centering
  \centerline{\includegraphics[width=\linewidth]{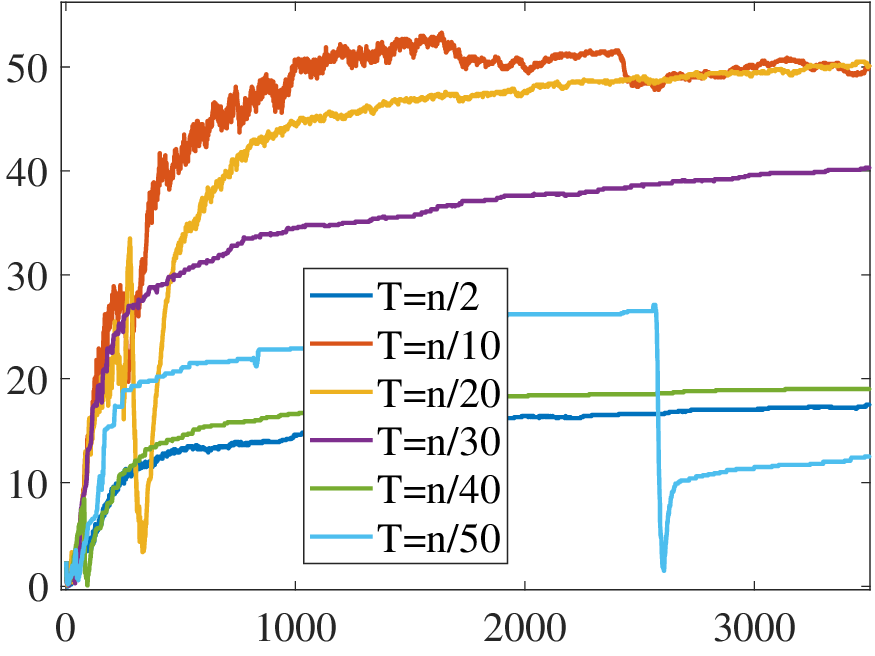}}
  \centerline{(b) \name-V2 - sparsity} \medskip
\end{minipage}

\begin{minipage}[b]{.49\linewidth}
  \centering
  \centerline{\includegraphics[width=\linewidth]{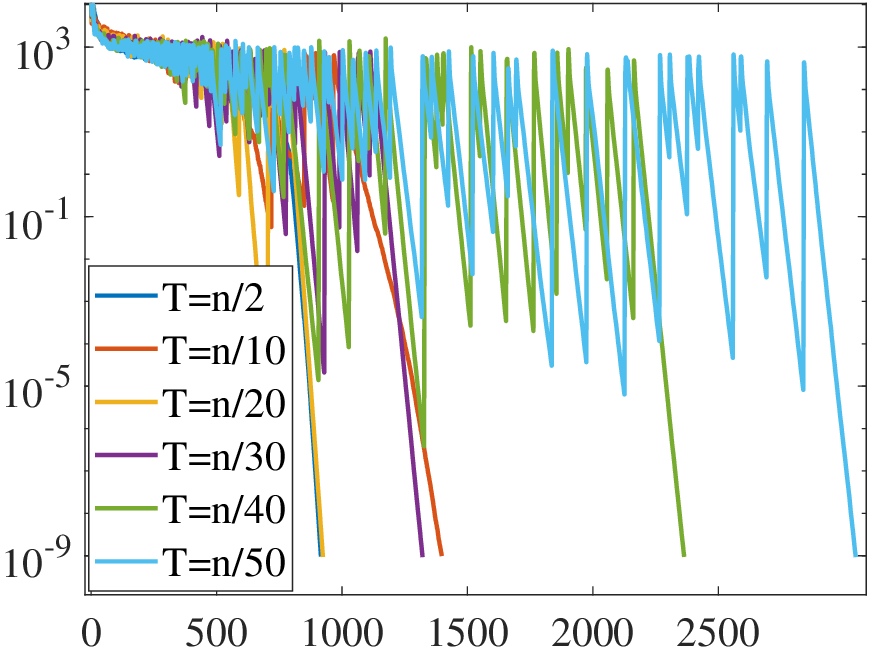}}
  \centerline{(a) \name-V3 - gradient norm}\medskip
\end{minipage}
\hfill
\begin{minipage}[b]{0.49\linewidth}
  \centering
  \centerline{\includegraphics[width=\linewidth]{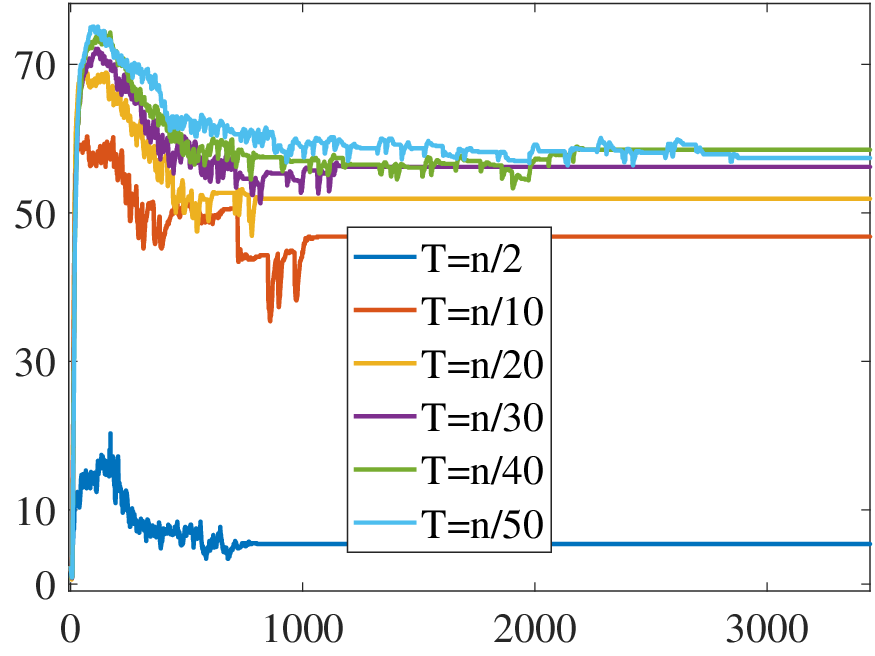}}
  \centerline{(b) \name-V3 - sparsity} \medskip
\end{minipage}
\caption{Norm of the gradient on the left and percentage of parameters below $\delta=10^{-3}$ on the right along the iterations for \name-V2 and \name-V3 and different target numbers of relevant parameters $T$. (Random least-squares A1)}
\label{Fig:change_T}
\end{figure}

As can be expected, the choice of the target number of relevant parameters $T$ does influence the behaviour of the four versions of \name. As we now show for problem A1 in Figure~\ref{Fig:change_T}, this effect may vary from version to version. Indeed, asking for a small $T$ does not necessarily result in a large final number of parameters with small magnitude, as is clear for the (unscaled) \name-V2.  Fortunately, the behaviour of the (scaled) \name-V3 is typically more consistent. For problem A1 and in general, choosing a smaller $T$ then results in a larger number of solution components whose absolute value is below the sparsity threshold, but, as can be expected, at the price of slower convergence (see the two bottom panels of Figure~\ref{Fig:change_T}). 

\subsubsection{\name, FW and Adagrad}

We now turn to results obtained when running the four versions of \name\ along FW and the standard Adagrad,
using $\tau_1=50$ for FW1, $\tau_2=100$, and $\beta=0.001$ for FW2. Table \ref{tab:rand_prun} in Apppendix~\ref{app:RLS} reports a comparison between the four versions of \name, FW1, FW2, and Adagrad, and displays, for each random matrix and each algorithm, the pruning quality measure $\rho$ defined in \req{qual_meas}, averaged over 20 runs, for 5 different percentages $\sigma$ of pruned components. In summary, these results show
that, when the percentage of pruned components $\sigma$ is below 50$\%$, \name-V3 is the most reliable method in five of the six problems considered, followed by \name-V4. On the contrary,
\name-V1 is extremely reliable when the percentage of pruned parameters is below $30\%$, but its performance degrades rapidly as this percentage increases. For very aggressive pruning, that is for $\sigma$ around $70\%$, FW1, FW2 and \name-V4 exhibit the best results, making those methods particularly suitable for applications where sparsity is to be preferred to high accuracy. 

\begin{figure}[ht]
\centering
\begin{minipage}[b]{.49\linewidth}
  \centering
  \centerline{\includegraphics[width=\linewidth]{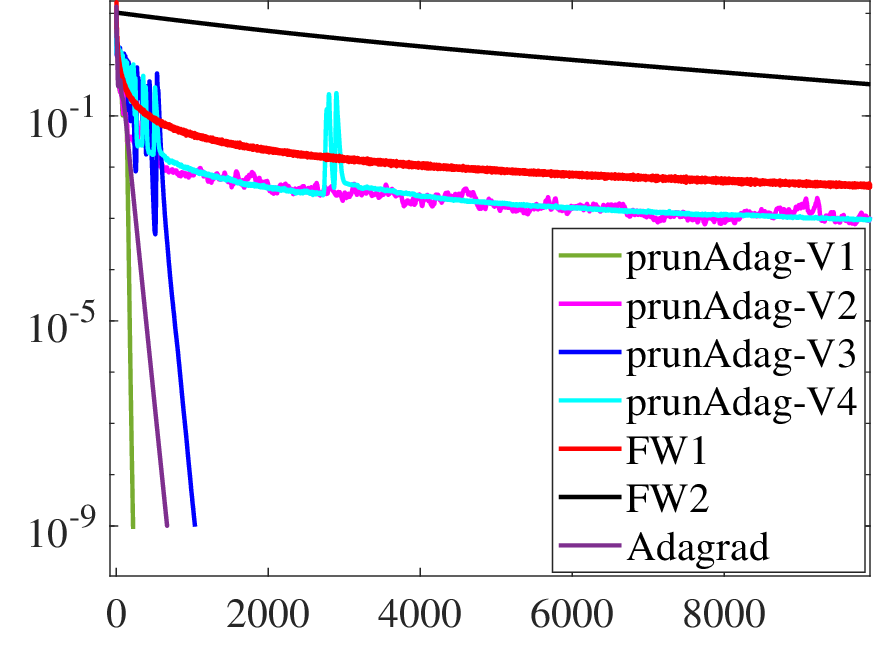}}
  \centerline{(a) gradient norm}\medskip
\end{minipage}
\hfill
\begin{minipage}[b]{0.49\linewidth}
  \centering
  \centerline{\includegraphics[width=\linewidth]{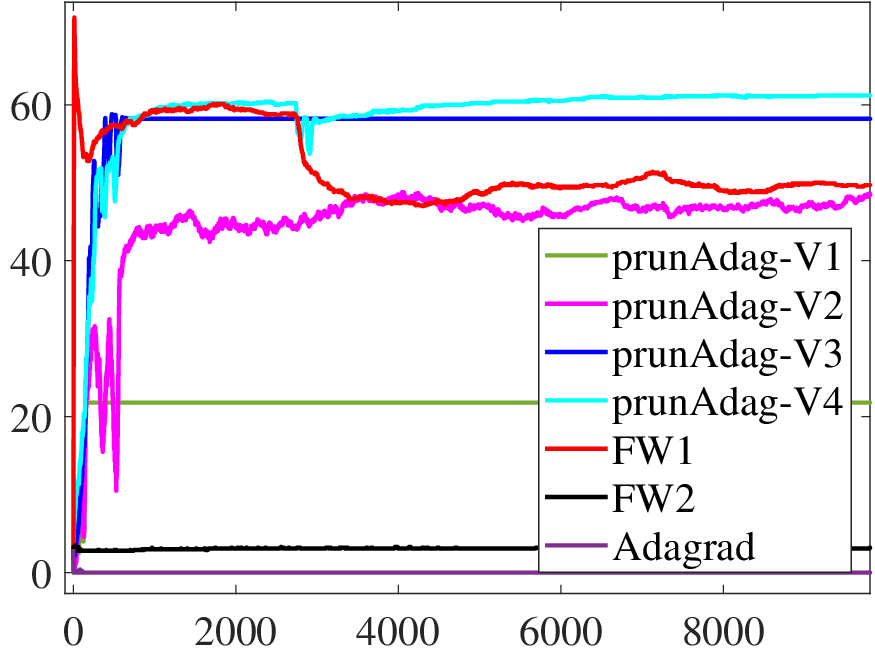}}
  \centerline{(b) sparsity} \medskip
\end{minipage}

\begin{minipage}[b]{0.49\linewidth}
  \centering
  \centerline{\includegraphics[width=\linewidth]{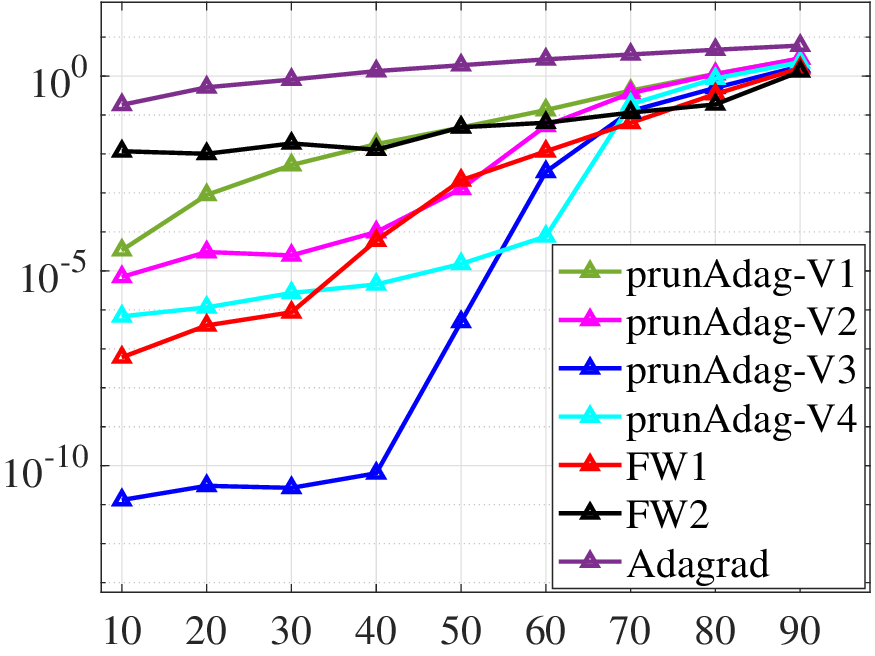}}
  \centerline{(c) $\omega$ measure} \medskip
\end{minipage}
\begin{minipage}[b]{0.49\linewidth}
  \centering
  \centerline{\includegraphics[width=\linewidth]{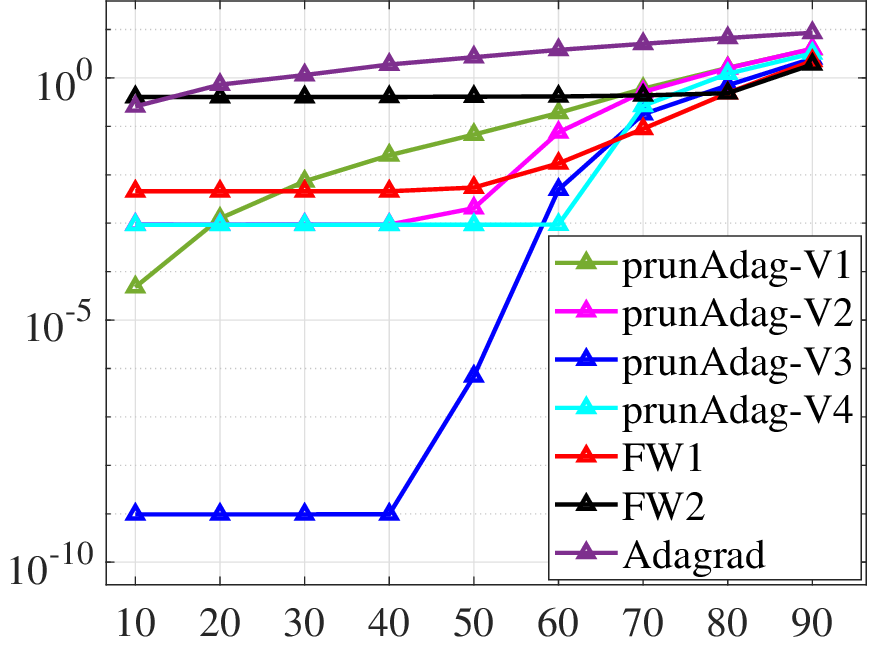}}
  \centerline{(d) $\rho$ measure} \medskip
\end{minipage}

\caption{On top, (a) gradient norm and (b) percentage of components below a fixed threshold $\delta=10^{-3}$ along iterations; at the bottom, (c) error measure $\omega$ and (d) error measure $\rho$ for different percentages of pruned components after the optimization. (Random least-squares A2)}
\label{Fig:rand_grad_spar}
\end{figure}

We illustrate these results by graphically detailing, in Figure \ref{Fig:rand_grad_spar}, the (typical) results obtained for problem A2.

This figure shows that version  \name-V1 is the fastest algorithm, performing comparably to the original Adagrad algorithm. However, it is the less robust to pruning among all \name\ versions. \name-V3 satisfies the stopping criterion on the gradient norm while being by far the best choice in terms of robustness to pruning up to $50\%$ of sparsity. 
Despite their poor performance for lower percentages of pruned components, FW1 and FW2 remain a valid alternative for very aggressive pruning, even though their convergence is the slowest among the algorithms considered. However, one should remember that FW is quite sensitive to the choice of its parameters and those used here have been tuned once for all the least-square problems considered (faster convergence can sometimes be achieved by further problem-by-problem tuning).

\subsection{SPARCO problems}

The aim of sparse signal recovery is finding a sparse representation of an observed noisy signal $b$ as a linear combination of some redundant dictionary $A$. Typically, the dictionary is a wide matrix with more columns than rows, consisting of various bases such as wavelet, discrete cosine, and Fourier. The SPARCO library \cite{van2007sparco} as supplied by S2MPJ \cite{GratToin24} includes examples of these problems for different dictionaries. Given a sparse vector $x^*$, the observation is generated as $b=Ax^*+r$, where $r $ is additive noise vector of appropriate dimension and $A$ is a fixed dictionary. The aim is to recover a robust solution by solving the related under-determined least-squares and using the \name\ algorithm instead of enhancing sparsity using a regularizing term\footnote{Since our focus is robustness to pruning, we deliberately ignore the robustness-to-noise issue that might occur if we solve the non-regularized least-squares. However, as long as the assumptions AS.1, AS.2, and AS.3 are satisfied, any regularization term can be added to the objective function.} 
For these tests, we set the FW parameter to $\tau_1=100$, $\tau_2=100$ and $\beta=0.001$. The complete results for different percentages $\sigma$ of pruned components are given in Table \ref{tab:sparco_prun_tab} in Appendix~\ref{app:sparco}. We observe that \name-V3 is the most robust algorithm in four out of five considered examples when the percentage of pruned components is below $50\%$. As for random-least squares, FW1 and \name-V3 are the most reliable for a percentage of pruned components $\sigma$ exceeding $50\%$, thus they represent a better choice for very aggressive pruning.

\begin{figure}[ht]
\centering 
\begin{minipage}[b]{.49\linewidth}
  \centering
  \centerline{\includegraphics[width=\linewidth]{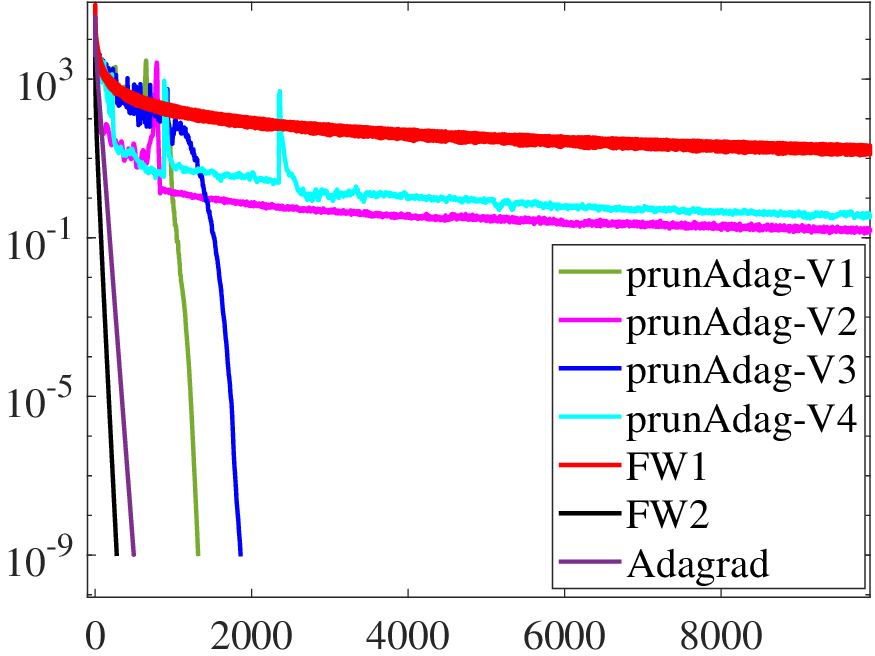}}
  \centerline{(a) gradient norm}\medskip
\end{minipage}
\hfill
\begin{minipage}[b]{0.49\linewidth}
  \centering
  \centerline{\includegraphics[width=\linewidth]{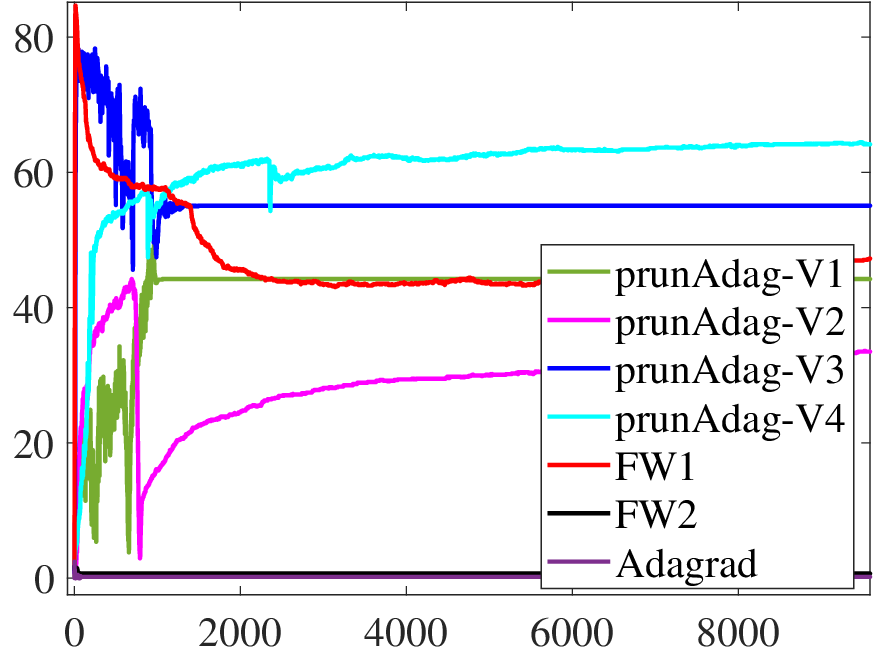}}
  \centerline{(b) sparsity} \medskip
\end{minipage}

\begin{minipage}[b]{0.49\linewidth}
  \centering
  \centerline{\includegraphics[width=\linewidth]{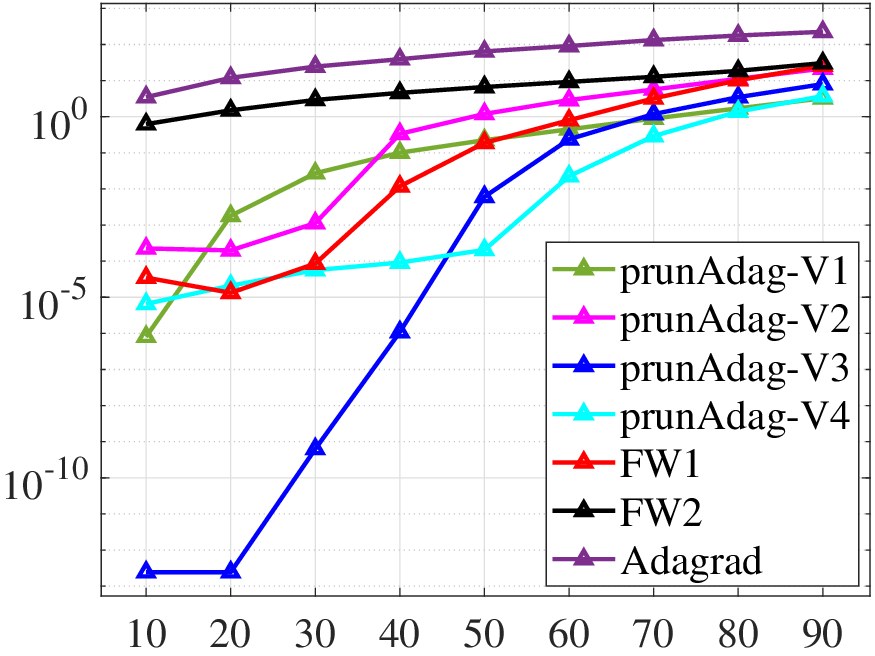}}
  \centerline{(c) $\omega$ measure} \medskip
\end{minipage}
\begin{minipage}[b]{0.49\linewidth}
  \centering
  \centerline{\includegraphics[width=\linewidth]{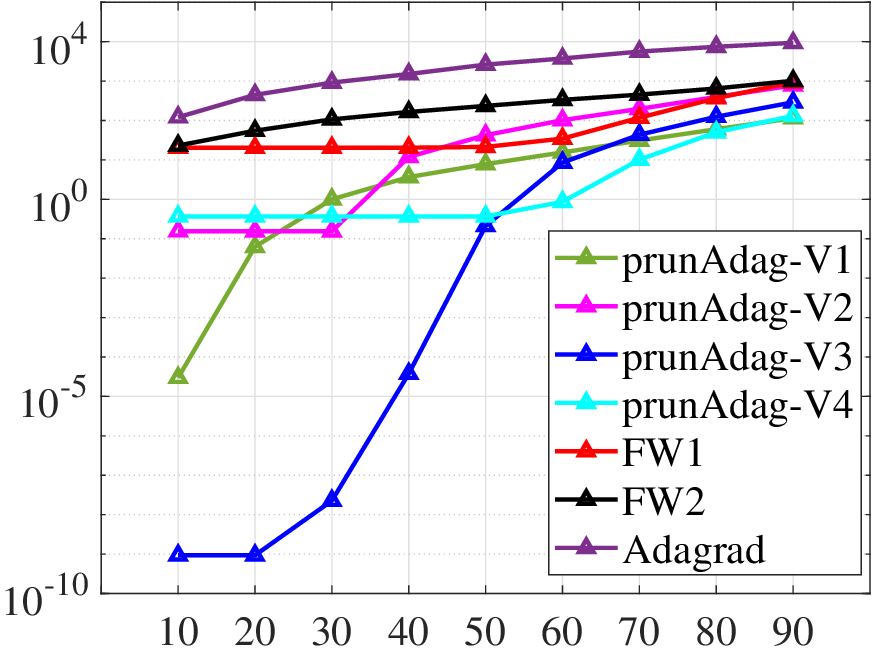}}
  \centerline{(d) $\rho$ measure} \medskip
\end{minipage}

\caption{On top, (a) gradient norm and (b) percentage of components below the threshold $\delta=10^{-3}$ along iterations; at the bottom, (c) error measure $\omega$ and (d) error measure $\rho$ for different percentages of pruned components $\sigma$ after the optimization. (Sparco 11)}
\label{Fig:sparco_grad}
\end{figure}

Figure \ref{Fig:sparco_grad} shows that FW2 is the fastest algorithm to converge, however, it does not exhibit strong robustness properties, showing a similar behaviour to Adagrad algorithm. 
By contrast, \name-V3 reaches the tolerance set for the norm of the gradient and it has the lowest value of the error measure $\rho$ up to $50\%$ of pruned components and the lowest value for error measure $\omega$ up to $40\%$ of pruned components. 
All versions of \name\ exhibit more robust performance than Adagrad algorithm. 

\clearpage

\subsection{Sparse coding in dictionary learning}

Let $Y$ be a given dataset, the aim of dictionary learning is to find a dictionary $D$ and a sparse coefficient matrix $X$, such that $Y \approx DX$. This problem is frequently solved by alternating optimization and we focus on the so-called sparse coding step, that is, given $D$ we aim at finding a sparse $X$ such that $Y \approx DX$. Given $k>0$, for each element of the dataset $y$ (column of $Y$) and positive integer $c$, the standard formulation of the sparse coding step is the following
\begin{equation}
\min_x \rVert y - Dx \lVert^2 \qquad \text{ such that } \rVert x\lVert_0 \leq c,
 \label{sparse_coding}
 \end{equation}
where $\ell_0$ denotes the zero-norm of a vector, defined as the number of its nonzero entries. We test our framework by addressing problem (\ref{sparse_coding}) removing the explicit constraint and using \name\ to find a possibly dense solution that is robust to pruning. This approach has the advantage of allowing \textit{a posteriori} pruning with different sparsity levels.

\begin{figure}[ht]
\begin{minipage}[h]{.24\linewidth}
  \centering
  \centerline{\includegraphics[width=\linewidth]{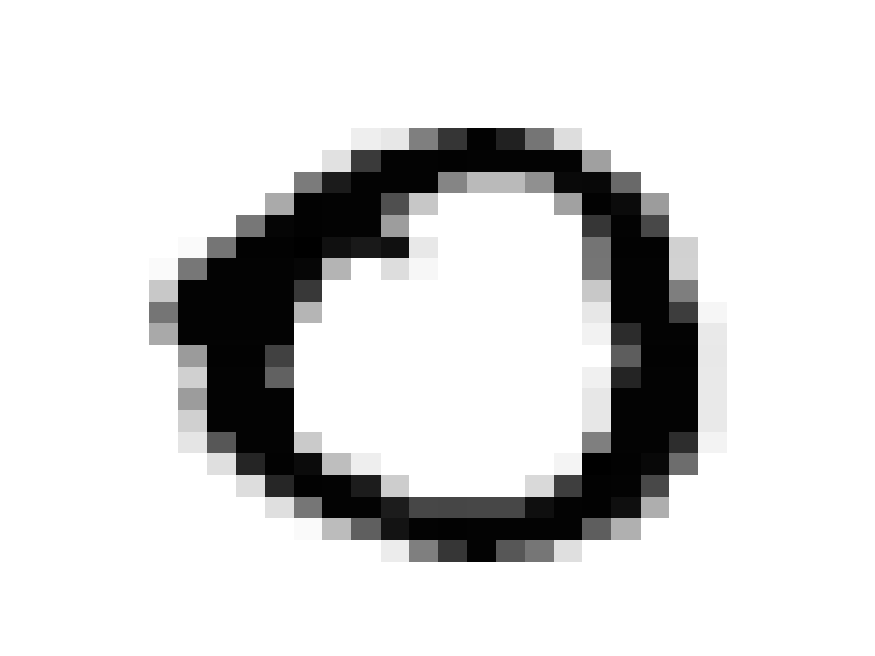}}
  \centerline{Real}
\end{minipage}
\begin{minipage}[h]{0.24\linewidth}
  \centering
  \centerline{\includegraphics[width=\linewidth]{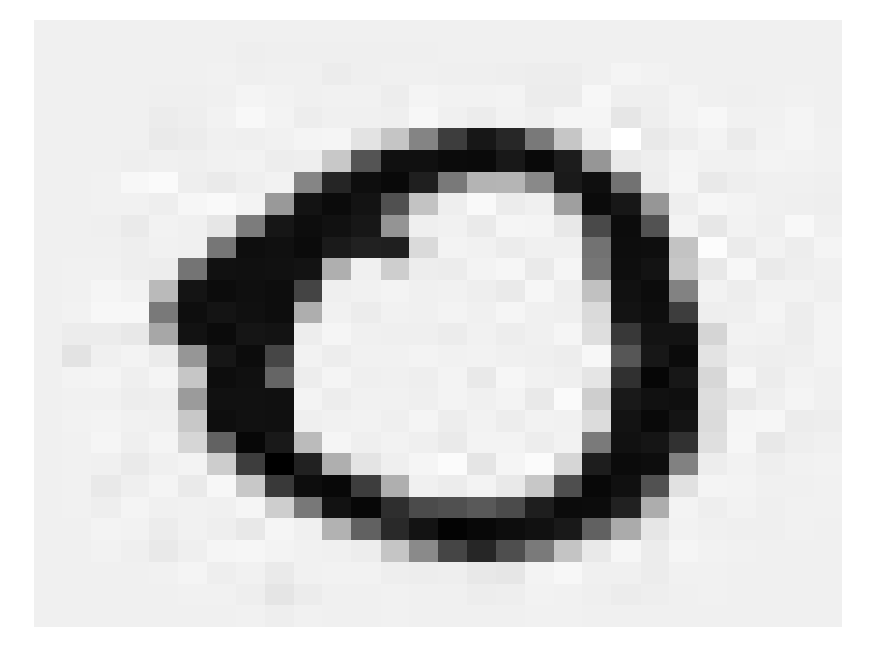}}
  \centerline{0$\%$ sparse}
\end{minipage}
\begin{minipage}[h]{0.24\linewidth}
  \centering
  \centerline{\includegraphics[width=\linewidth]{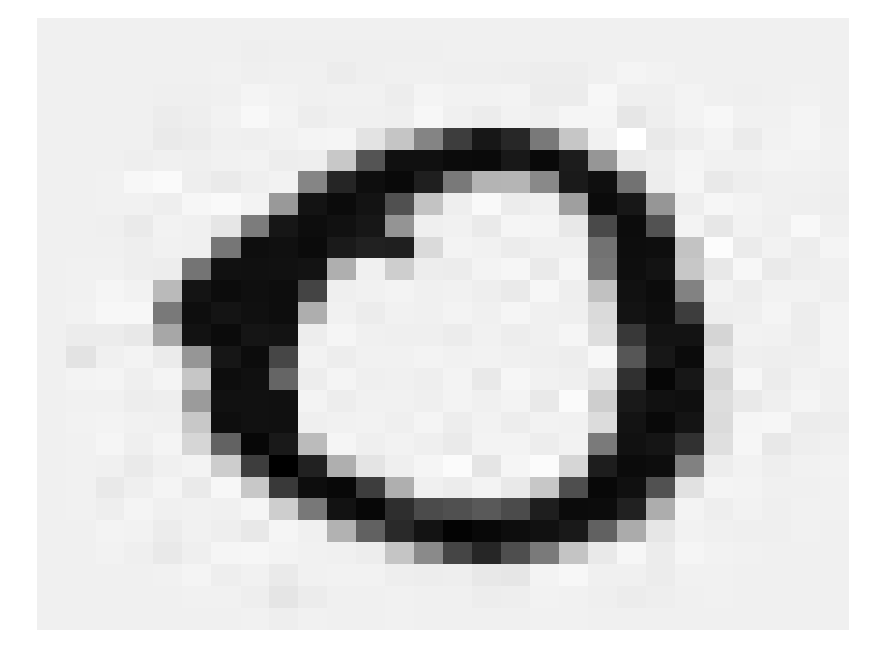}}
  \centerline{20$\%$ sparse}
\end{minipage}
\begin{minipage}[h]{0.24\linewidth}
  \centering
  \centerline{\includegraphics[width=\linewidth]{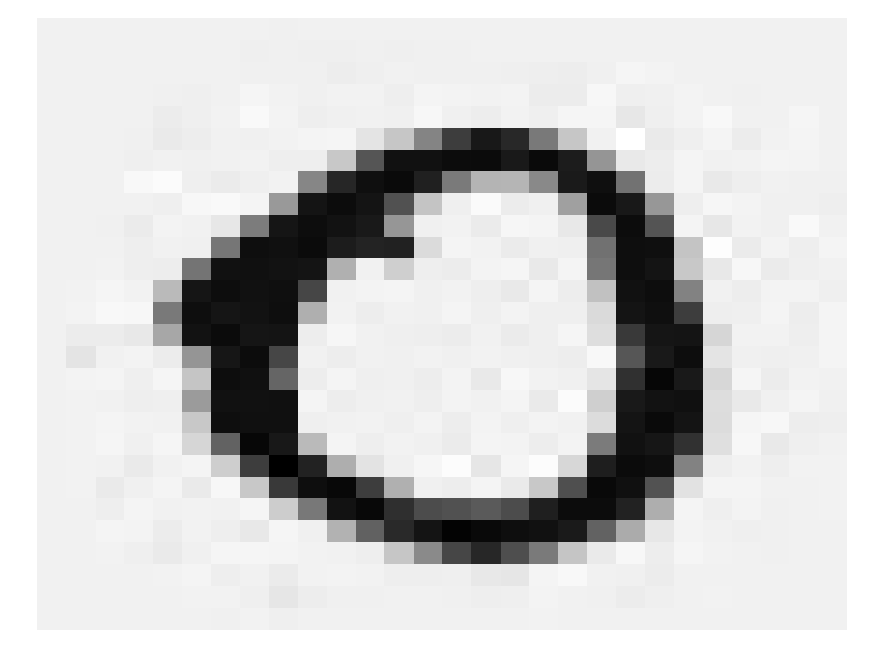}}
  \centerline{40$\%$ sparse}
\end{minipage}
\begin{minipage}[h]{.24\linewidth}
  \centering
  \centerline{\includegraphics[width=\linewidth]{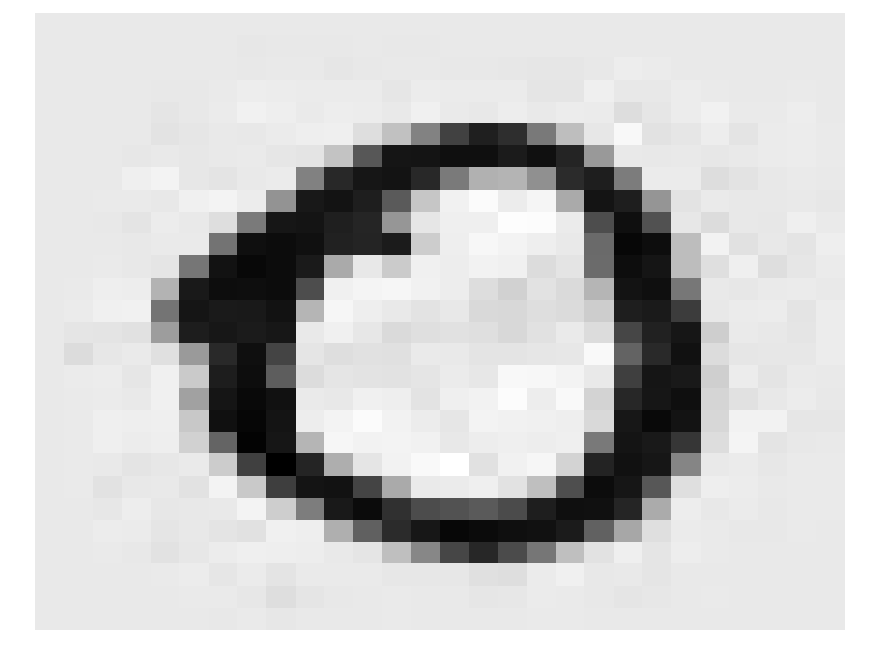}}
  \centerline{60$\%$ sparse}
\end{minipage}
\begin{minipage}[h]{0.24\linewidth}
  \centering
  \centerline{\includegraphics[width=\linewidth]{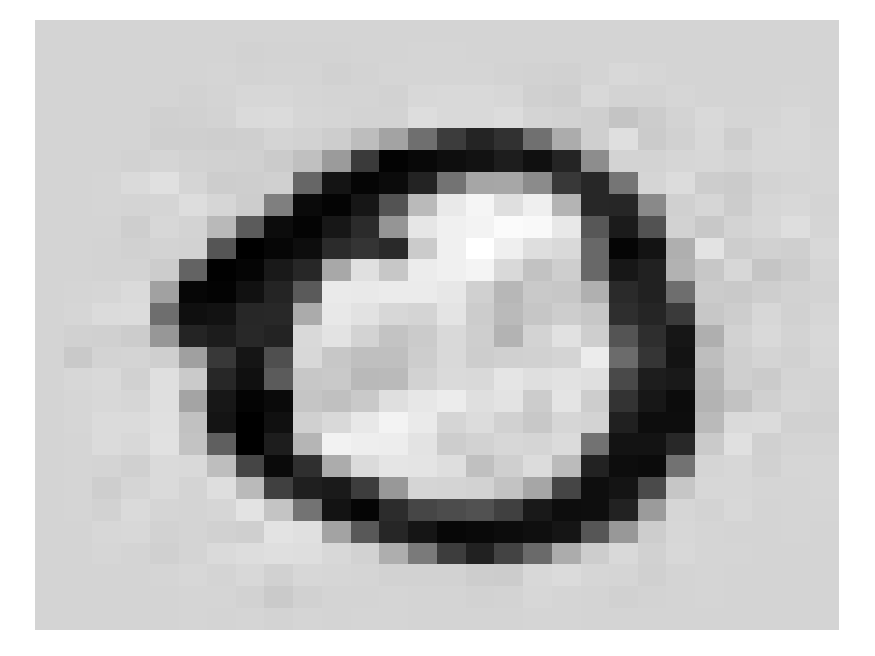}}
  \centerline{70$\%$ sparse}
\end{minipage}
\begin{minipage}[h]{0.24\linewidth}
  \centering
  \centerline{\includegraphics[width=\linewidth]{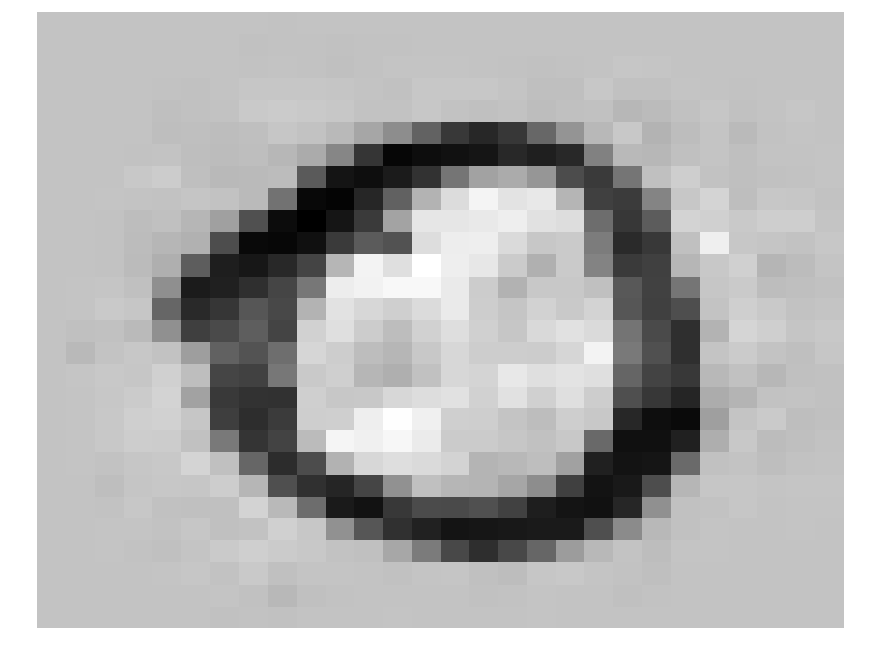}}
    \centerline{80$\%$ sparse}
\end{minipage}
\begin{minipage}[h]{0.24\linewidth}
  \centering
  \centerline{\includegraphics[width=\linewidth]{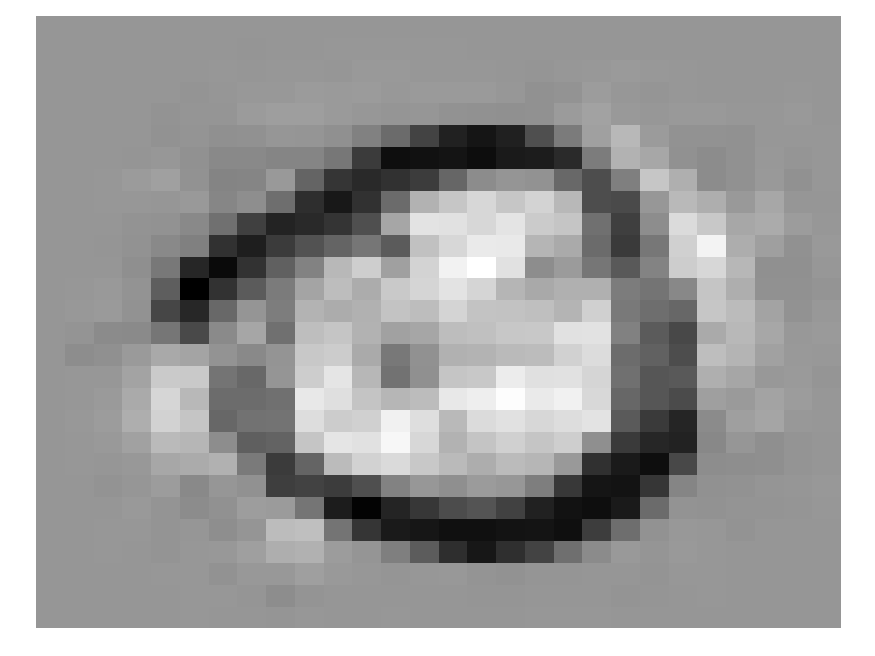}}
  \centerline{90$\%$ sparse}
  \end{minipage}
  \medskip
  
 \centerline{(a) \name-V1} \medskip
 
\begin{minipage}[h]{.24\linewidth}
  \centering
  \centerline{\includegraphics[width=\linewidth]{./img_real.eps}}
  \centerline{Real}
\end{minipage}
\begin{minipage}[h]{0.24\linewidth}
  \centering
  \centerline{\includegraphics[width=\linewidth]{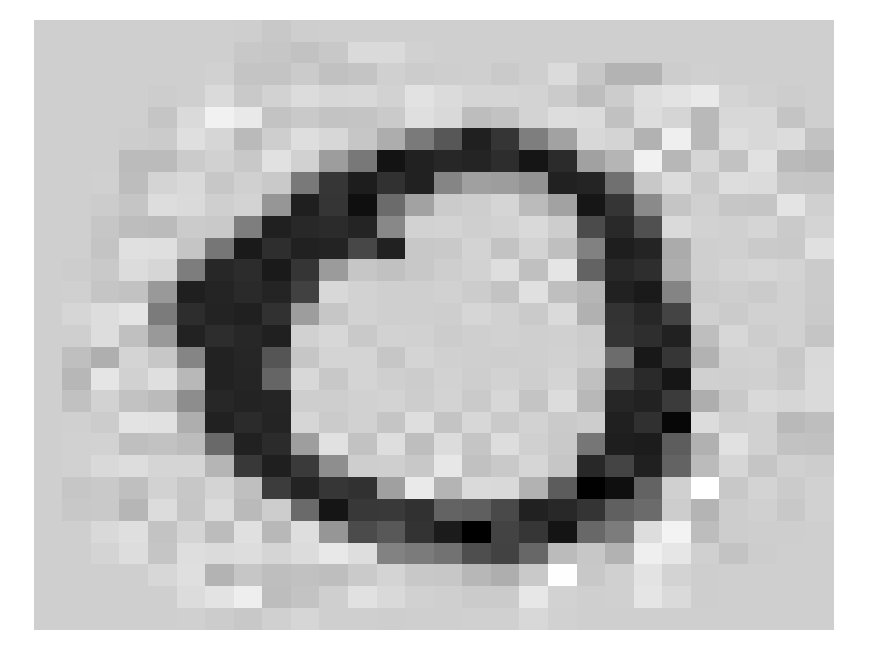}}
  \centerline{ 0$\%$ sparse}
\end{minipage}
\begin{minipage}[h]{0.24\linewidth}
  \centering
  \centerline{\includegraphics[width=\linewidth]{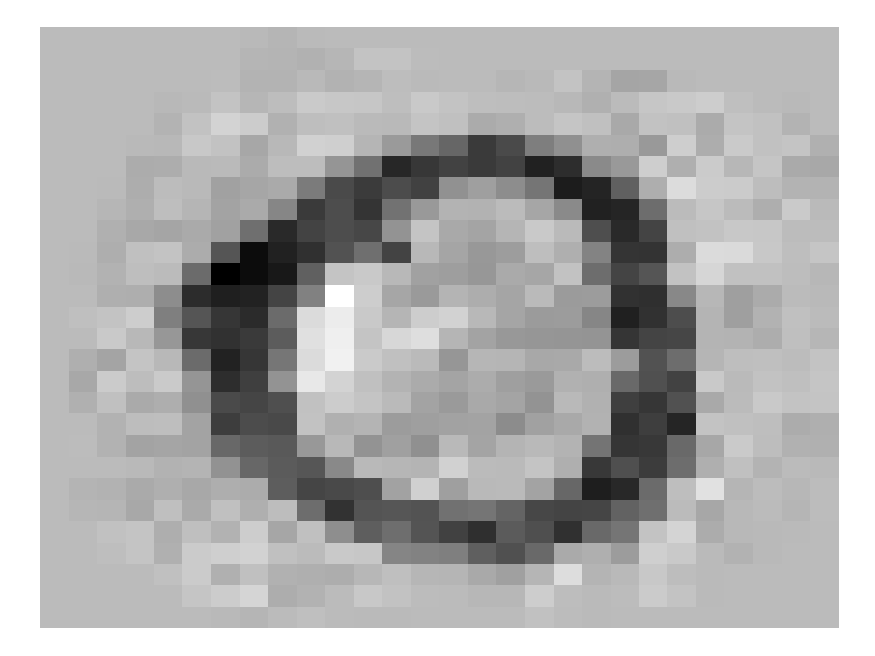}}
  \centerline{20$\%$ sparse}
\end{minipage}
\begin{minipage}[h]{0.24\linewidth}
  \centering
  \centerline{\includegraphics[width=\linewidth]{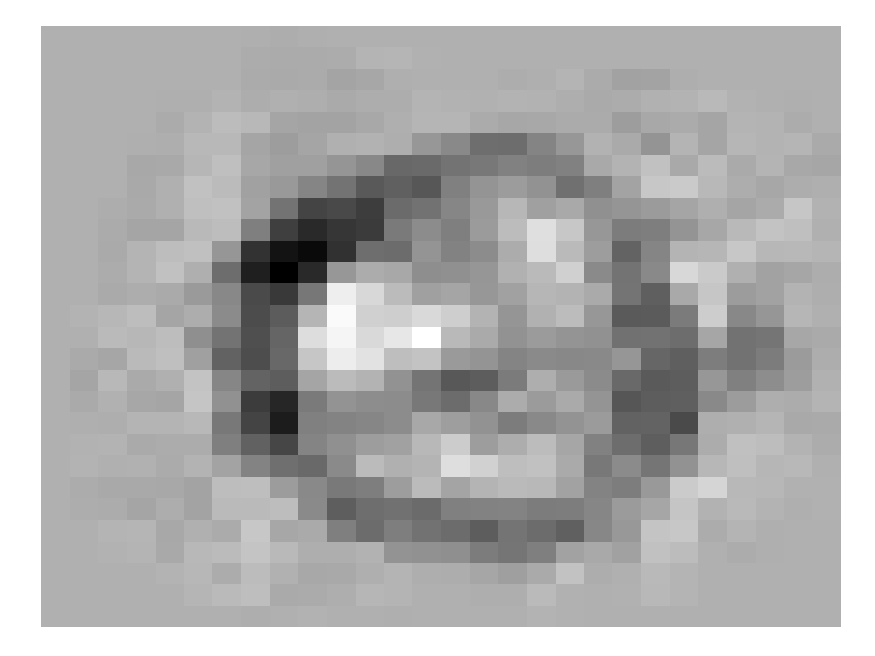}}
  \centerline{40$\%$ sparse}
\end{minipage}
\begin{minipage}[h]{.24\linewidth}
  \centering
  \centerline{\includegraphics[width=\linewidth]{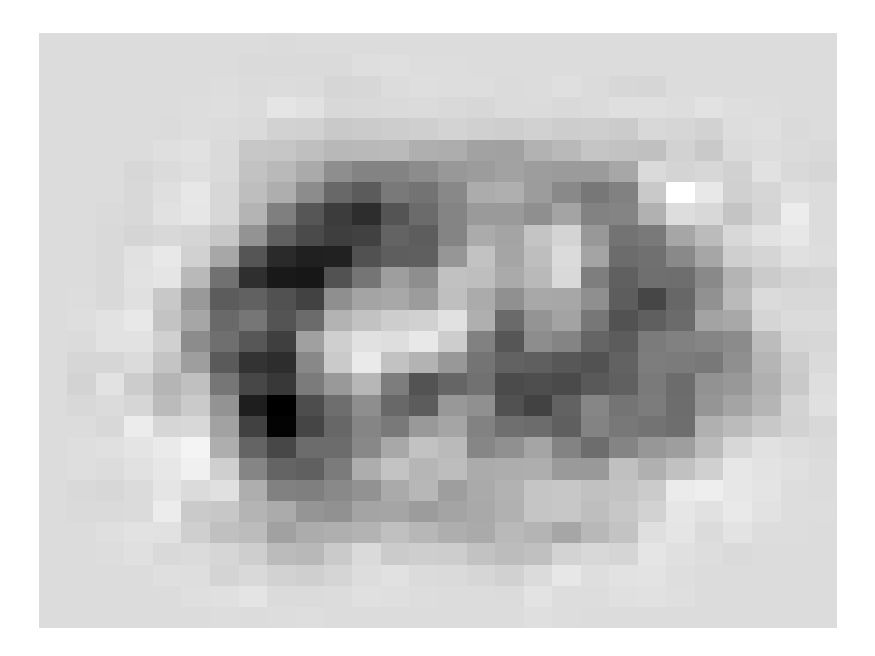}}
  \centerline{60$\%$ sparse}
\end{minipage}
\begin{minipage}[h]{0.24\linewidth}
  \centering
  \centerline{\includegraphics[width=\linewidth]{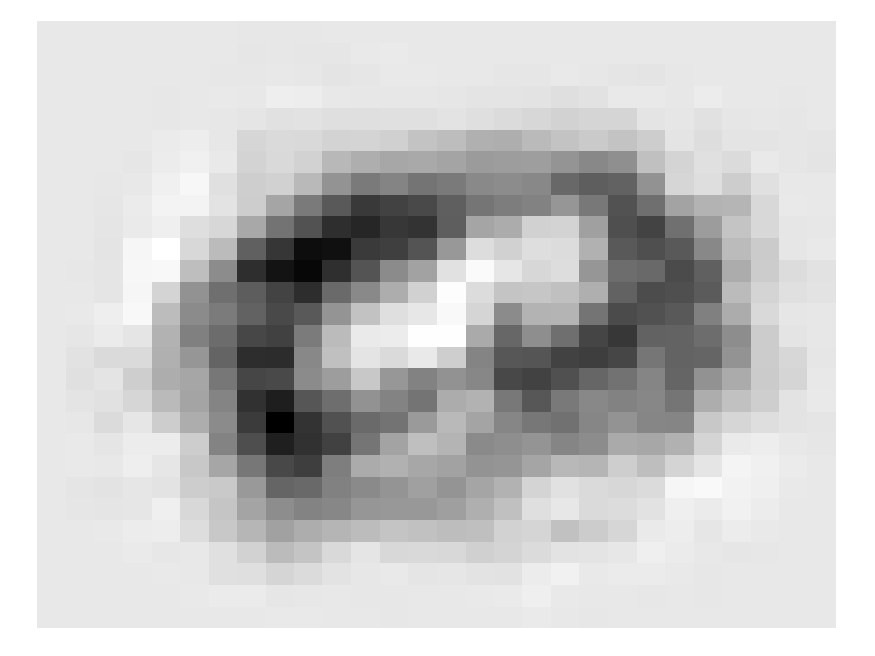}}
  \centerline{70$\%$ sparse}
\end{minipage}
\begin{minipage}[h]{0.24\linewidth}
  \centering
  \centerline{\includegraphics[width=\linewidth]{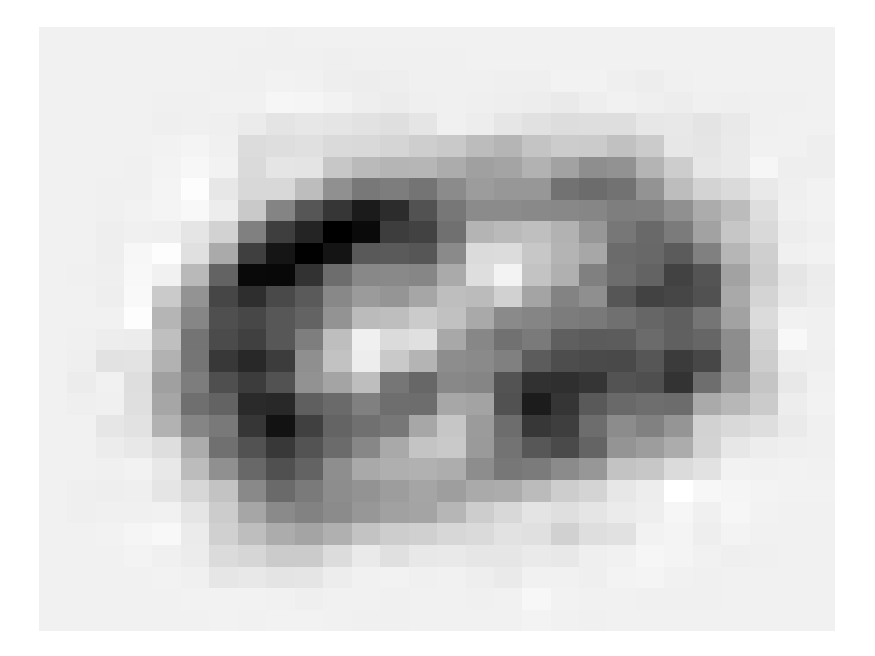}}
    \centerline{80$\%$ sparse}
\end{minipage}
\begin{minipage}[h]{0.24\linewidth}
  \centering
  \centerline{\includegraphics[width=\linewidth]{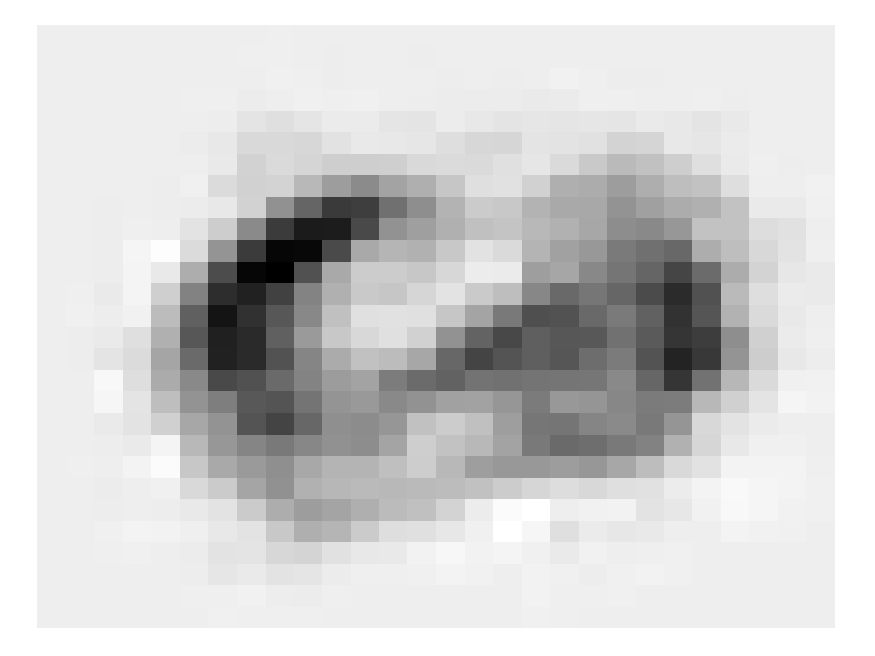}}
  \centerline{90$\%$ sparse}
\end{minipage}
 \medskip
 \centerline{(b) Adagrad}  
\caption{MNIST sparse coding step in dictionary learning. Visual representations of pruned solutions of our best algorithm \name-V1 in dictionary learning application, for different percentage $\sigma$ of pruned parameters in the first two rows (a). Adagrad solution for different levels of pruned components in the last two rows (b).}
\label{Fig:method_comp_dict}
\end{figure}

\begin{figure}[h]
 \centering
\begin{minipage}[b]{.49\linewidth}
  \centerline{\includegraphics[width=\linewidth]{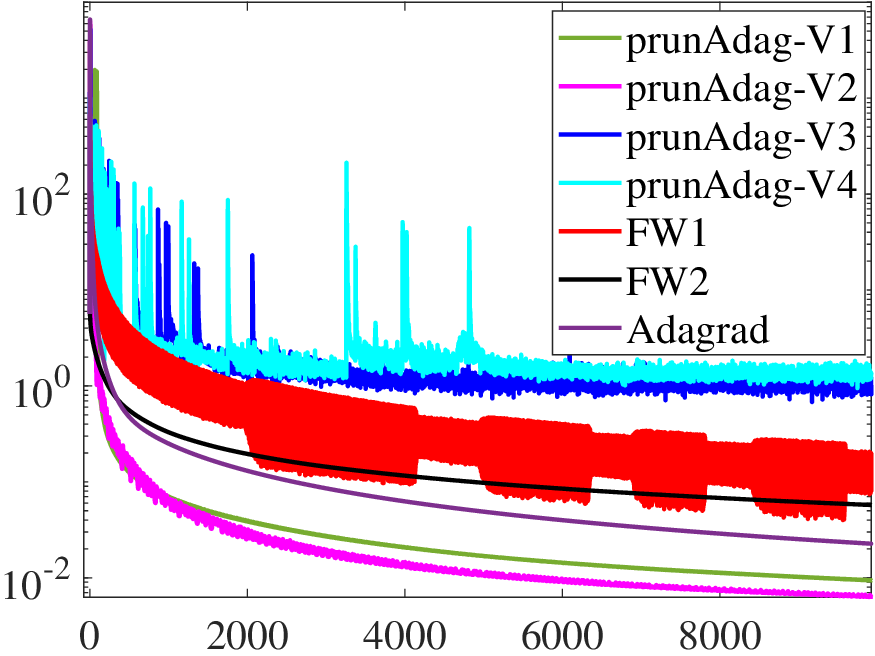}}
  \centerline{(a) gradient norm} \medskip
\end{minipage}
\begin{minipage}[b]{0.49\linewidth}
  \centering
  \centerline{\includegraphics[width=\linewidth]{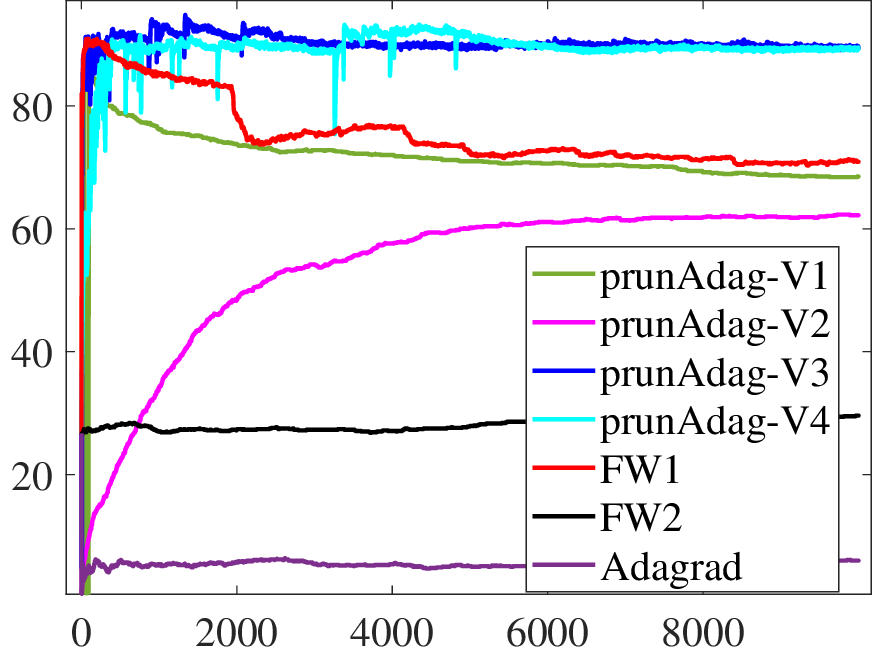}}
  \centerline{(b) sparsity} \medskip
\end{minipage}

\begin{minipage}[b]{0.49\linewidth}
  \centering
  \centerline{\includegraphics[width=\linewidth]{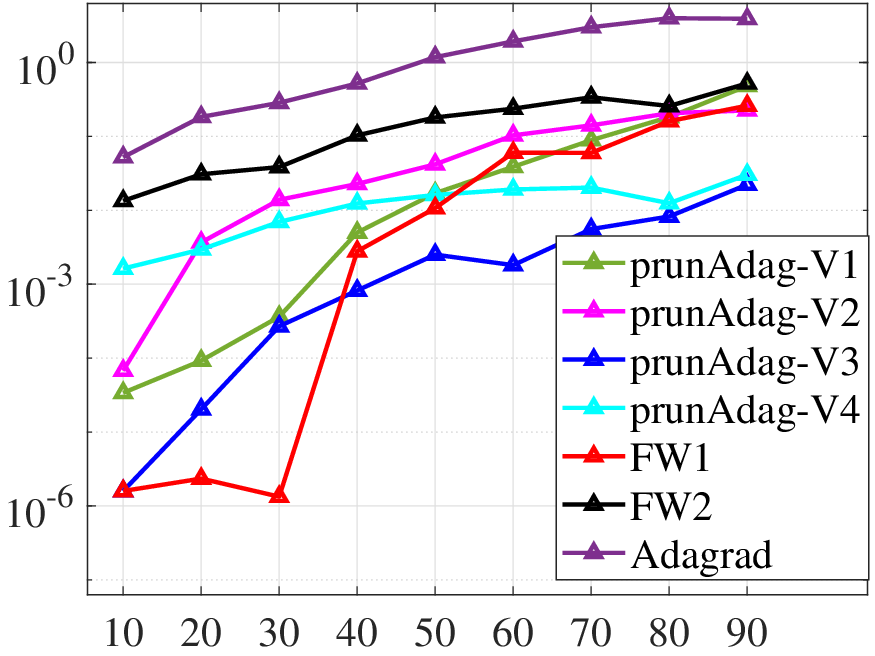}}
  \centerline{(c) $\omega$ measure} \medskip
\end{minipage}
\begin{minipage}[b]{0.49\linewidth}
  \centering
  \centerline{\includegraphics[width=\linewidth]{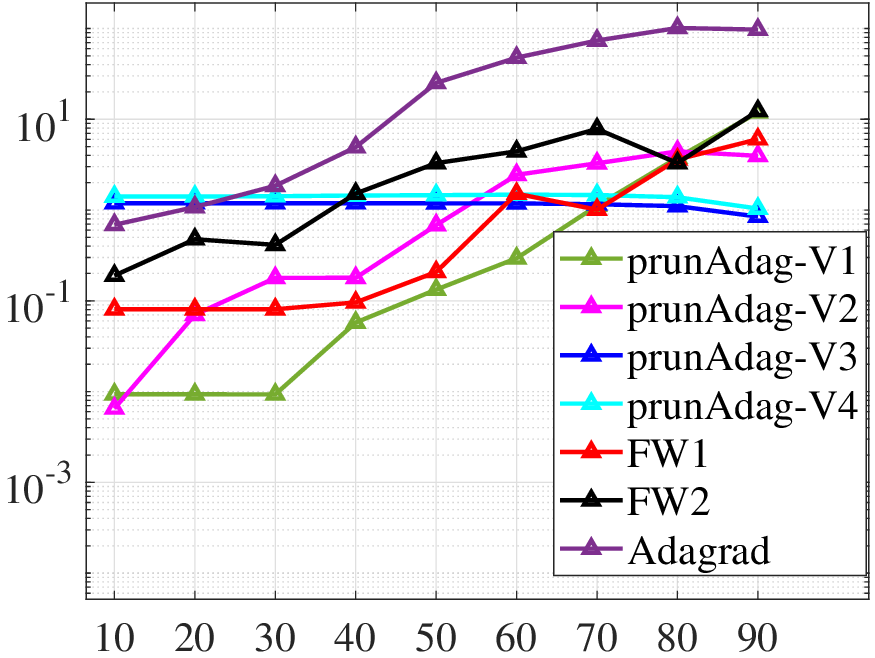}}
  \centerline{(d) $\rho$ measure} \medskip
\end{minipage}

\caption{On top, (a) gradient norm decrease and (b) percentage of parameters below the threshold $\delta=10^{-2}$  along iterations; at the bottom, (c) error measure $\omega$ and (d) error measure  $\rho$ for pruned components percentage $\sigma$ from $10\%$ to $90\%$ on the $x$-axis (MNIST).}
\label{Fig:dict_analysis}
\end{figure}

In our tests, we consider a subset of the MNIST data set~\cite{lecun1989optimal} and we generated $D$ in problem \req{sparse_coding} by using KSVD \cite{ksvd}, the state-of-the-art solver for solving the dictionary learning problems\footnote{We used the Matlab implementation KSVD-Box v13 of K-SVD available at http://www.cs.technion.ac.il/$\sim$ronrubin/software.html with default parameters.}. The data set $Y$ has dimension $784\times4000$, the dictionary $D$ has 784 rows and 1000 columns and $c$ in
\req{sparse_coding} is chosen equal to 100. We set the FW parameters to $\tau_1=10$, $\tau_2=20$, and $\beta=0.001$.
In Figure \ref{Fig:method_comp_dict} we illustrate the reconstruction of an instance $y \approx D\Bar{x}$, where $\Bar{x}$ is the pruned solution obtained by \name-V1 algorithm for increasing percentages of pruned parameters. The results show that the solution begins to degrade when more than $40\%$ of parameters are pruned. Figure \ref{Fig:dict_analysis} confirms the visual intuition since the error measure $\rho$ remains below $10^{-1}$ for \name-V1 algorithm. Moreover, Figure \ref{Fig:dict_analysis} on the left highlights that \name-V1 is the most robust algorithm up to $70\%$ of pruned components, while \name-V3 and \name-V4, despite their poor global accuracy, remain the best choice for aggressive pruning with more than $80\%$ of pruned parameters.

\clearpage

\subsection{Binary classification}

\begin{figure}[ht]
\begin{minipage}[b]{.49\linewidth}
  \centering
  \centerline{\includegraphics[width=\linewidth]{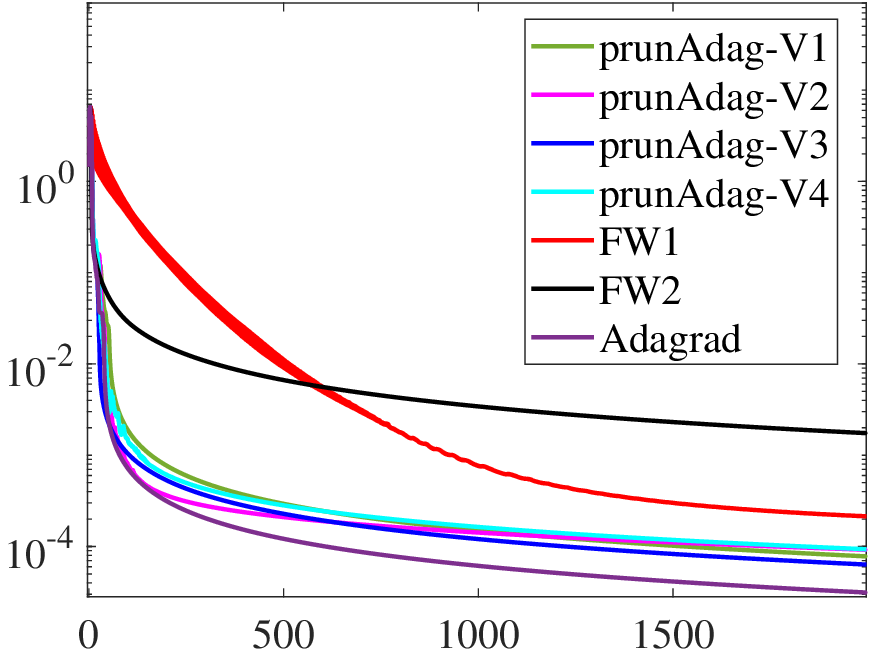}}
  \centerline{(a) gradient norm}\medskip
\end{minipage}
\hfill
\begin{minipage}[b]{0.49\linewidth}
  \centering
  \centerline{\includegraphics[width=\linewidth]{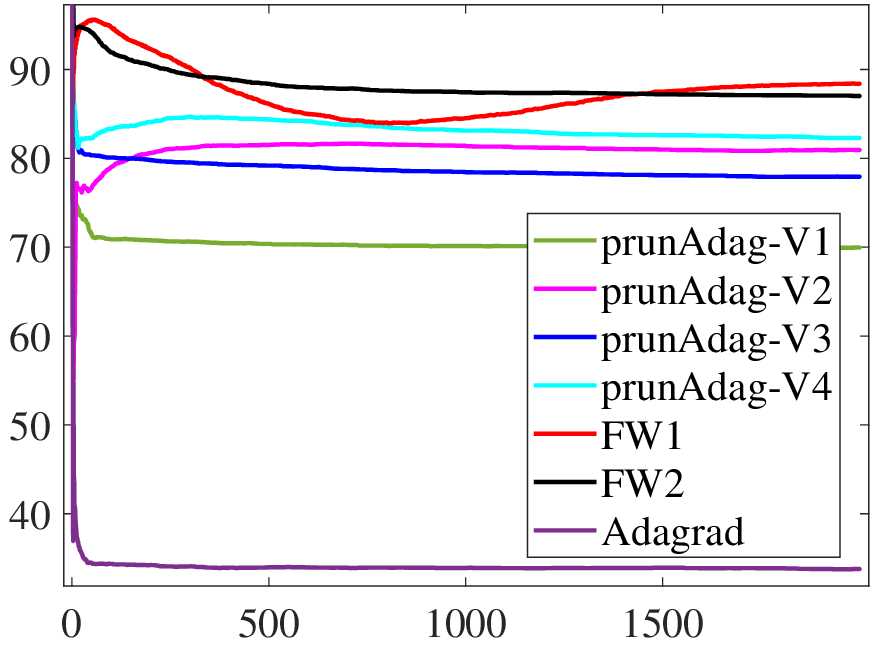}}
  \centerline{(b) sparsity} \medskip
\end{minipage}
\centering
\begin{minipage}[b]{0.49\linewidth}
  \centering
  \centerline{\includegraphics[width=\linewidth]{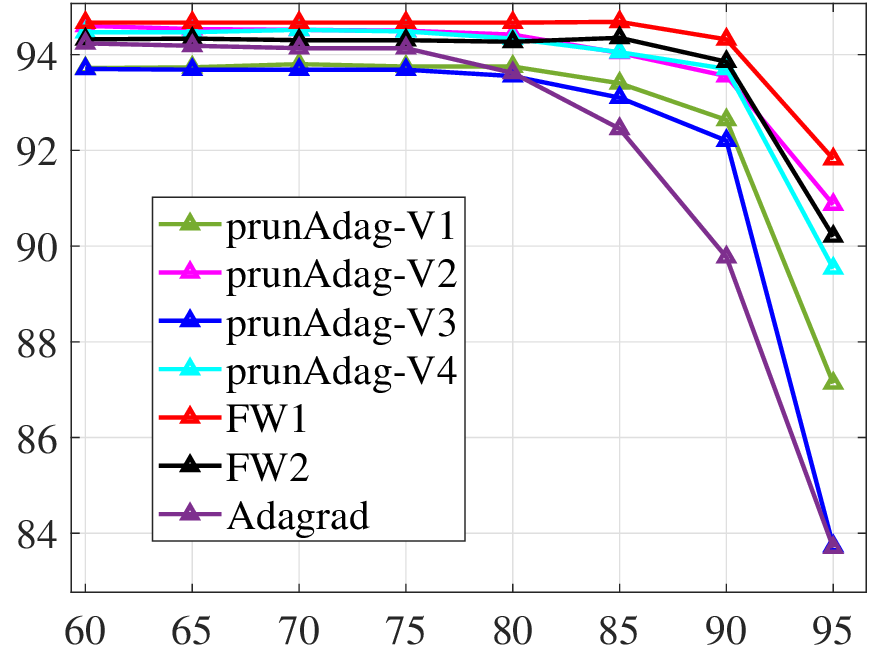}}
  \centerline{(c) classification accuracy} \medskip
\end{minipage}
\caption{On top, (a) gradient norm decrease and (b) percentage of parameters below the threshold $\delta=10^{-1}$ along iterations; 
at the bottom, (c) average percentage of correctly classified samples as a function of sparsity $\sigma$ from $60\%$ and $95\%$ (GISETTE).
}
\label{Fig:grad_gisette_fw}
\end{figure}

Finally, we test our method on the averaged logistic loss for binary classification. We assume that a labeled training set $\{a_i,b_i\}$ with $a_i \in \mathbb{R}^n$ and $b_i \in \{0,1\}$ for $i=1,\dots,N$ is available, where $b_i$ classifies each sample into two distinct classes. The averaged logistic loss over all samples is neither linear nor convex, and it is defined as 
\begin{equation}
f(x)=\frac{1}{N} \sum_{i=1}^N \log(1+e^{-b_i a_i^T x}).
\label{logistic}
\end{equation}
If the number of features in the data set is large, we expect that some may be redundant or irrelevant in the classification process; thus, pruning can be used to achieve a sparse solution that does not consistently degrade the classification performance. We use \name\ to minimize the function in (\ref{logistic}) on the training set and to promote convergence towards a solution $x$ in which the largest components correspond to relevant features. Then, we prune the parameters to achieve different levels of sparsity $\sigma$ and evaluate the prediction on the testing set using the pruned solution. For this experiment, we select small-size data sets\footnote{We randomly selected 1000 samples for MNIST and A9A.} that have at least 100 features, that are 
MNIST\footnote{Classification between even and odd numbers.}~\cite{lecun1989optimal}, GISETTE~\cite{machrep}, 
REGEDO~\cite{chang2011libsvm}, 
A9A~\cite{machrep}, 
and MOLECULE~\cite{machrep}.
We split the training and testing set following a ratio 70:30. The data is normalized using min-max normalization and each algorithm is randomly initialized and stopped after 2000 iterations. We set the FW parameters to $\tau_1=10,$ $\tau_2=100$, and $\beta=0.5$. The results are collected in Table \ref{tab:bin_class_tab} in Appendix~\ref{app:bin_class}, where we show the average test accuracy and the percentage of components pruned for each algorithm. In Figure \ref{Fig:grad_gisette_fw} we analyze the performance of the algorithms on the GISETTE data set.  

Table~\ref{tab:bin_class_tab} in Appendix~\ref{app:bin_class} shows that all four versions of \name\ yield a consistent reduction in the number of parameters of the model over 80$\%$ for all problems, without affecting the classification performance, and are significantly more robust to pruning than Adagrad. We are not able to identify one of the four versions that outperforms the others; however, we can observe that \name-V4 tolerates the largest number of pruned components. 
 In Figure \ref{Fig:grad_gisette_fw} (c) we compare Adagrad, FW and \name\ in terms of average robustness to pruning parameters after training on the GISETTE data set from 20 independent starting points. 
 The figure clearly shows that all four versions of \name are more robust to pruning than Adagrad, as small components correspond to irrelevant components of the model. Indeed, all the versions show a stable accuracy for a $70\%$ sparse solution and that of \name-V2 and \name-V4 are not affected significantly by pruning $90\%$ of solution's components.  One also notes the excellent performance of FW on this example and level of sparsity. It seems that, even if the norm of the gradient at the final iterate is still larger after the training phase than is the case for \name\ (see Figure~\ref{Fig:grad_gisette_fw} (a)), it is sufficient to produce good model predictions.

 \clearpage
 
\section{Conclusions} \label{sec:concl}

We have proposed a new first-order OFFO method, named \name, intended for applications where pruning of the variables/parameters is desirable. The new ``pruning-aware" algorithm uses a new strategy to classify parameters at each iteration of \name\ algorithm into ``optimisable" and ``decreasable", instead of ``relevant" and ``irrelevant" as suggested in~\cite{ding2019global}, and extending the concept introduced in~\cite{zimmer2022compression}, where the optimization is performed on the components related to the largest partial derivatives.
It also features a new framework to update parameters in these two classes separately, based on an Adagrad-like step for the first and on an adaptive trust-region approach to decrease the magnitude of the variables in the second. We proved the convergence of the method to first-order stationary points with global rate $O(\log(k)/\sqrt{k+1})$. Finally, we conducted numerical experiments on several real-world applications, such as sparse signal recovery, dictionary learning, and binary classification. These experiments suggest that the new approach (and its \name-V3\ version in particular) has a clear practical potential.

While we have, in this paper, focused on the ``deterministic case" where the gradient values are 
computed exactly, the ``stochastic case" where the gradient may be contaminated by random noise 
(such as sampling) is also clearly of interest, and the object of current research. Also of 
interest is the inclusion of momentum in \name\ or a similar algorithm.

{\footnotesize
\subsection*{\footnotesize Acknowledgment}

Ph.~T. gratefully acknowledges the friendly support of DIEF (University of Florence) during his visit in the fall of 2024.

\subsection*{\footnotesize Ethics approval and consent to participate}
Not applicable.
\subsection*{\footnotesize Consent for publication}
Not applicable.
\subsection*{\footnotesize Funding}
{\footnotesize
The work of M.P. and G.S. was partially supported by INdAM-GNCS under the
project CUP E53C23001670001.  The research of M.P.  was  partially granted by the Italian Ministry of University and Research (MUR) through the PRIN 2022 ``MOLE: Manifold constrained Optimization and LEarning'',  code: 2022ZK5ME7 MUR D.D. financing decree n. 20428 of November 6th, 2024 (CUP B53C24006410006), and by PNRR - Missione 4 Istruzione e Ricerca - Componente C2 Investimento 1.1, Fondo per il Programma Nazionale di Ricerca e Progetti di Rilevante Interesse Nazionale (PRIN) funded by the European Commission under the NextGeneration EU programme, project ``Advanced optimization METhods for automated central veIn Sign detection in multiple sclerosis from magneTic resonAnce imaging (AMETISTA)'',  code: P2022J9SNP,
MUR D.D. financing decree n. 1379 of 1st September 2023 (CUP E53D23017980001).
G.S. acknowledges the support by the European Union (ERC consolidator, eLinoR, no 101085607).
}
\subsection*{\footnotesize Availability of data and materials}
The data that support the findings of this study are available from the corresponding author upon request.
\subsection*{\footnotesize Competing interests}
The authors declare that they have no competing interests.
\subsection*{\footnotesize Authors' contributions}
All authors contributed equally to the writing of this article. All authors reviewed the manuscript.

\appendix
\clearpage
\section{Detailed numerical results}
\subsection{Under-determined least-squares}\label{app:RLS}

\vspace*{-3mm}
\begin{table}[ht]
\centering
\begin{tabular}{|c|c|c|c|c|c|c|c|c|}
\hline
\hspace{-0.2cm} Problem \hspace{-0.2cm} & $\sigma$& \hspace{-0.3cm} V1 \hspace{-0.3cm} &V2 &V3 &V4 &FW1 &FW2&Adagrad \\
\hline
\multirow{5}{*}{A1} & $10\%$&  $2.2e^{-5}$& $1.39$   & \boldmath{$9.4e^{-10}$} & $0.88$ & $6.34$ & $8.92$& $670$ \\
\cline{2-9}
                       & $20\%$ &  0.01 & 1.39  & \boldmath{$9.7e^{-10}$} & 0.88 & 6.34 & 25.67& $1.7e^{3}$\\
\cline{2-9}
                       &$30\%$  &  0.21& 1.38  & \boldmath{$5.2e^{-4}$} & 0.88 & 6.34 & 47.2& $3.0e^{3}$\\
\cline{2-9}
                       & $40\%$&  1.69& 1.36  & \boldmath{$0.17$} & 0.88 & 6.43 & 79.1& $4.4e^{3}$\\
\cline{2-9}
                       & $50\%$&  12.9& 1.35  & 9.73 & \textbf{0.89} & 7.52 & 109& $6.0e^{3}$\\
\hline
\multirow{5}{*}{A2} & $20\%$&  $8.7e^{-4}$ & $8.2e^{-4}$  & \boldmath{$5.9e^{-6}$} & $1.1e^{-3}$ & $4.5e^{-3}$ & $0.39$& 0.61 \\
\cline{2-9}
                       & $40\%$ &  0.02 & $8.2e^{-4}$  & \boldmath{$6.0e^{-6}$} & $1.1e^{-3}$ & $5.4e^{-3}$ & $0.40$& 1.72\\
\cline{2-9}
                       &$60\%$  &  0.20 &  0.12 & $5.3e^{-3}$ & \boldmath{$2.2e^{-3}$} & $0.03$& 0.42&3.57\\
\cline{2-9}
                       & $70\%$& 0.63 & $0.57$  & 0.15 & 0.21 & \textbf{0.14} & $0.44$& 4.90\\
\cline{2-9}
                       & $80\%$&  1.62 & 1.67  & 0.67 & 1.06 & 0.60 & \textbf{0.49}& 6.43\\
\hline
\multirow{5}{*}{A3} & $20\%$&  $4.1e^{-3}$ & 0.02  & \boldmath{$9.7e^{-10}$} & $9.7e^{-3}$ & 0.06 & $0.09$& 6.29 \\
\cline{2-9}
                      & $30\%$&  0.12 & 0.02  & \boldmath{$9.7e^{-10}$} & $9.7e^{-3}$ & 0.06 & $0.44$& 29.7 \\
\cline{2-9}
                       & $40\%$ &  0.38 & 0.02  & \boldmath{$3.4e^{-5}$} & $9.7e^{-3}$ & 0.06& 0.70& 44.1\\
\cline{2-9}
                       & $50\%$ & 1.19 & 0.01  & 0.03 & \boldmath{$9.8e^{-3}$} & 0.08 & 1.00& 60.1\\
\cline{2-9}
                       & $60\%$&  3.57 & 0.33  & 0.50 & 0.23 & \textbf{0.22} & 1.38& 76.7\\
\hline
\multirow{5}{*}{A4} & $20\%$&  $6.1e^{-4}$ & $8.8e^{-4}$  & \boldmath{$4.2e^{-8}$} & $9.7e^{-4}$ & $4.5e^{-3}$ & 0.38& 0.65 \\
\cline{2-9}
                       & $30\%$ &  $4.1e^{-3}$ &  $8.8e^{-4}$  & \boldmath{$4.2e^{-8}$} & $9.7e^{-4}$ & $4.6e^{-3}$ & 0.38& 1.13\\
\cline{2-9}
                       &$50\%$  &  0.06 &  $6.9e^{-3}$ &\boldmath{$7.5e^{-6}$} & $9.7e^{-4}$ & 0.01 & 0.39&2.65\\
\cline{2-9}
                       & $70\%$&  0.65 &  0.57 &\textbf{0.13} & 0.22& 0.14 & 0.42&5.12\\
\cline{2-9}
                       & $90\%$&  3.98 & 4.10  & 2.51 & 3.34 & 2.53 & \textbf{1.91}& 8.43\\
\hline
\multirow{5}{*}{A5} & $10\%$&  $3.7e^{-4}$ & 1.58  & \boldmath{$9.5e^{-10}$} & 0.88 & 6.34& 8.96& 648 \\
\cline{2-9}
                       & $20\%$&  0.03 & 1.58  & \boldmath{$9.5e^{-10}$} & 0.88 & 6.35& 25.4& $1.7e^{3}$ \\
\cline{2-9}
                       & $30\%$ &  0.34 & 1.57  & \boldmath{$6.5e^{-4}$} & 0.88 & 6.35& 47.6& $3.0e^{3}$\\
\cline{2-9}
                      &$40\%$  &  2.41 &  1.54 & \textbf{0.12} & 0.88 & 6.47& 77.0&$4.4e^{3}$\\
\cline{2-9}
                       & $50\%$& 14.6 & 1.49  & 7.64 & \textbf{0.89} & 8.21 & 113& $6.0e^{3}$\\

\hline
\multirow{5}{*}{A6} & $10\%$&  $2.2e^{-3}$ & \boldmath{$4.3e^{-4}$}  & $3.0e^{-3}$ & $3.6e^{-4}$ & $6.4e^{-3}$& 1.40& 0.20 \\
\cline{2-9}
                       & $30\%$&  0.11 & 0.13  & 0.29 & 0.10 & \textbf{0.02}& 1.40& 1.34 \\
\cline{2-9}
                       & $50\%$ &  0.64 & 0.80  & 1.01 & \textbf{0.73} & 0.35& 1.41& 3.15\\
\cline{2-9}
                      &$60\%$  &  1.23 &  1.47 & 1.60 & 1.40 & \textbf{0.91}& 1.42&4.2\\
\cline{2-9}
                       & $70\%$& 2.2 & 2.5  & 2.4 & 2.4 & 1.89 & \textbf{1.44}& 5.33\\

\hline
\end{tabular}
\caption{Random least-squares: the error measure $\rho$ for different percentages of pruned components $\sigma$.}
\label{tab:rand_prun}
\end{table}
\clearpage

\subsection{SPARCO problems}\label{app:sparco}
\vspace*{-3mm}
\begin{table}[h]
\centering
\begin{tabular}{|c|c|c|c|c|c|c|c|c|}
\hline
Problem& $\sigma$& V1 &V2 &V3 &V4 &FW1 &FW2&Adagrad \\
\hline
\multirow{5}{*}{S3} & $10\%$&  $4.2e^{-5}$&$3e^{-3}$   &\boldmath{$9.6e^{-10}$} & $1.6e^{-3}$ & 0.02 & 0.41& 1.52 \\
\cline{2-9}
                       & $30\%$ &  $8.1e^{-3}$  &$3.1e^{-3}$  & \boldmath{$1.2e^{-4}$} & $1.6e^{-3}$& 0.02 & 0.41& 7.02\\
\cline{2-9}
                       &$50\%$  &  0.08& $0.05$ & 0.12 & \boldmath{$2.7e^{-3}$} & 0.12 & 0.41& 13.5\\
\cline{2-9}
                       & $60\%$&  0.25& 0.08  & 0.55 & \boldmath{$6.1e^{-3}$} & 0.12 & 0.41& 17.2\\
\cline{2-9}
                        & $70\%$&  6.21& 6.69  & 6.29 & 5.56 & \textbf{0.31} & 0.41& 21.9\\
\hline
\multirow{5}{*}{S5} & $20\%$&  0.02 & 0.03  & \boldmath{$3.6e^{-5}$} & 0.02 & 0.11 & 0.01& 7.68\\
\cline{2-9}
                       & $30\%$ &  0.06 & 0.03  & \boldmath{$1.5e^{-3}$} & 0.02 & 0.11 & 0.22& 57.0\\
\cline{2-9}
                       &$40\%$  &  0.23 & \textbf{0.02} & 0.03 & \textbf{0.02} & 0.14& 0.42&88.3\\
\cline{2-9}
                       & $50\%$& 0.86 & 0.03  & 0.34 & \textbf{0.02} & 0.44 & 0.63& 125\\
\cline{2-9}
                       & $60\%$&  4.49 & 2.21  & 1.78 & 1.78 &1.83 & \textbf{0.87}& 164\\
\hline
\multirow{5}{*}{S7} & $10\%$& $2.5e^{-4}$ & \boldmath{$4.6e^{-6}$}  & $1.7e^{-5}$ & $6.0e^{-4}$ & 0.01 & 0.05& 0.08 \\
\cline{2-9}
                       & $30\%$ & $2.3e^{-3}$ & $1.4e^{-3}$  & \boldmath{$1.7e^{-5}$} & $6.0e^{-4}$ & 0.01 & 0.07& 0.43 \\
\cline{2-9}
                       &$50\%$  &  $7.0e^{-3}$ &  $7.0e^{-3}$  & \boldmath{$1.7e^{-5}$} & $6.0e^{-4}$ & 0.01 &  0.14& 0.95 \\
\cline{2-9}
                       & $70\%$&  0.02 & 0.02  & \boldmath{$2.0e^{-4}$} & $6.0e^{-4}$& 0.01 & 0.27& 2.0 \\
\cline{2-9}
                       & $90\%$&  0.03 & 0.05  & \boldmath{$3.3e^{-3}$} & \boldmath{$3.3e^{-3}$} & 0.01 & 0.50& 2.35 \\
\hline
\multirow{5}{*}{S9} & $10\%$&  \textbf{1.23} & 1.62  & 5.52 & 5.51 & 44.9 & 2.4 & 2.5 \\
\cline{2-9}
                       & $20\%$ & \textbf{1.23} & 1.62  & 5.53 & 5.51 & 44.9 & 9.56 & 9.46 \\
\cline{2-9}
                       &$30\%$ &  \textbf{1.23}  & 1.62  & 5.52 & 5.51 & 44.9 & 16.2 & 31.8 \\
\cline{2-9}
                       & $40\%$&  \textbf{1.23}  & 1.67  & 5.52 & 5.51 & 45.9 & 45.4 & 40.8 \\
\cline{2-9}
                       &$ 50\%$&  \textbf{1.23}  &  2.05 &5.46 & 5.54 & 33.6 & 75.7 &60.6\\
\hline
\multirow{5}{*}{S11} & $10\%$&  $3.3e^{-6}$ & 0.08 & \boldmath{$9.7e^{-10}$} & 0.50 & 19.3& 13.1& 165\\
\cline{2-9}
                       & $20\%$&  $4.2e^{-3}$ & 0.08 & \boldmath{$9.7e^{-10}$} & 0.50 & 19.3& 47.3& 489\\
\cline{2-9}
                       & $30\%$ &  0.13 & 0.07  & \boldmath{$2.5e^{-9}$} & 0.50 & 19.3& 84.3& 975 \\
\cline{2-9}
                       &$40\%$  &  1.19 &  1.01 & \boldmath{$3.9e^{-4}$} & 0.50 & 19.3&  152&$1.5e^{3}$\\
\cline{2-9}
                       & $50\%$&  5.6 &  12.1 & \textbf{0.37 }& 0.50 & 31.6&  332&$3.7e^{3}$\\
\hline
\end{tabular}
\caption{SPARCO problems: the error measure $\rho$ and different percentages of pruned components $\sigma$.}
\label{tab:sparco_prun_tab}
\end{table}

\vspace*{-5mm}

\clearpage
\subsection{Binary classification}\label{app:bin_class}
\vspace*{-3mm}
\begin{table}[h]
\centering
\begin{tabular}{|c|c|c|c|c|c|c|c|c|}
\hline
\hspace{-0.4cm} Problem \hspace{-0.4cm} & $\sigma$& V1 &V2 &V3 &V4 &FW1 &FW2&Adagrad \\
\hline
\multirow{5}{*}{GISETTE} & $75\%$&  $93.75\%$ &$94.50\%$   &$93.68\%$& $94.51\%$ & \boldmath{$94.66\%$} & $94.30\%$& $94.13\%$ \\
\cline{2-9}
                       & $80\%$ &  $93.75\%$ &$94.42\%$   &$93.55\%$& $94.33\%$ & \boldmath{$94.66\%$} & $94.26\%$& $93.61\%$ \\
\cline{2-9}
                       &$85\%$  &  $93.40\%$ &$94.03\%$   &$93.10\%$& $94.05\%$ & \boldmath{$94.68\%$} & $94.35\%$& $92.45\%$ \\
\cline{2-9}
                       & $90\%$&  $92.63\%$ &$93.55\%$   &$92.20\%$& $93.70\%$ & \boldmath{$94.31\%$} & $93.85\%$& $89.76\%$ \\
\cline{2-9}
                        & $95\%$&  $87.13\%$ &$90.86\%$   &$83.73\%$& $89.53\%$ & \boldmath{$91.81\%$} & $90.20\%$& $83.70\%$ \\
\hline
\multirow{5}{*}{MNIST} & $75\%$&  $81.90\%$ &$81.91\%$   &$80.60\%$& $80.46\%$ & \boldmath{$83.65\%$} & $81.85\%$& $80.15\%$ \\
\cline{2-9}
                       & $80\%$ &  $81.48\%$ &$81.53\%$   &$80.48\%$& $80.41\%$ & \boldmath{$83.33\%$} & $81.21\%$& $79.45\%$ \\
\cline{2-9}
                       &$85\%$  &  $80.43\%$ &$80.53\%$   &$80.35\%$& $80.12\%$ & \boldmath{$82.80\%$} & $80.83\%$& $77.32\%$ \\
\cline{2-9}
                       & $90\%$&  $78.02\%$ &$78.15\%$   &$78.83\%$& $78.85\%$ & \boldmath{$81.51\%$} & $79.02\%$& $75.68\%$ \\
\cline{2-9}
                        & $95\%$&  $72.26\%$ &$72.38\%$   &$74.11\%$& $73.92\%$ & \boldmath{$76.60\%$} & $72.48\%$& $66.06\%$ \\
\hline
\multirow{5}{*}{REGEDO} & $70\%$ & $96.36\%$  & $95.93\%$   & $96.36\%$  & $96.33\%$  & \boldmath{$97.33\%$}  & \boldmath{$97.33\%$}& $96.20\%$  \\
\cline{2-9}
                       & $75\%$ & $96.33\%$  & $95.93\%$   & $96.50\%$  & $96.67\%$  & \boldmath{$97.33\%$}  & \boldmath{$97.33\%$}& $95.86\%$  \\
\cline{2-9}
                      & $80\%$ & $96.20\%$  & $95.83\%$   & $96.23\%$  & $96.76\%$  & \boldmath{$97.33\%$}  & \boldmath{$97.33\%$}& $96.00\%$  \\
\cline{2-9}
                      & $85\%$ & $96.76\%$  & $96.10\%$   & $96.96\%$  & $96.60\%$  & $97.33\%$  & \boldmath{$97.46\%$}& $93.06\%$  \\
\cline{2-9}
                       & $90\%$ & $95.90\%$  & $95.67\%$   & $84.50\%$  & $96.40\%$  & $97.33\%$  & \boldmath{$98.60\%$}& $65.20\%$  \\
\hline
\multirow{5}{*}{A9A} & $65\%$& $81.65\%$  & $82.16\%$   & $81.87\%$  & $81.97\%$  & $83.07\%$  & \boldmath{$83.33\%$}& $76.93\%$  \\
\cline{2-9}
                       & $70\%$ & $81.33\%$  & $81.34\%$   & $81.53\%$  & $81.67\%$  & \boldmath{$83.37\%$}  & $83.10\%$& $72.78\%$  \\
\cline{2-9}
                       &$75\%$ & $80.68\%$  & $80.86\%$   & $80.77\%$  & $80.22\%$  & \boldmath{$83.40\%$}  & $83.03\%$& $71.25\%$  \\
\cline{2-9}
                       & $80\%$&  $79.53\%$  & $78.63\%$   & $79.55\%$  & $79.01\%$  & $82.70\%$  &  \boldmath{$82.85\%$}& $71.25\%$  \\
\cline{2-9}
                       &$ 85\%$& $75.83\%$  & $75.95\%$   & $76.20\%$  & $76.23\%$  & \boldmath{$82.36\%$}  & $81.60\%$& $72.76\%$  \\
\hline
\multirow{5}{*}{MOLECULE} &$ 65\%$& $79.30\%$  & $78.77\%$   & $78.70\%$  & $78.39\%$  & $78.53\%$  & \boldmath{$79.54\%$}& $65.00\%$  \\
\cline{2-9}
                       & $70\%$& $77.51\%$  & $78.11\%$   & $78.43\%$  & $78.18\%$  & $78.56\%$  & \boldmath{$78.95\%$}& $62.76\%$  \\
\cline{2-9}
                       & $75\%$& $76.11\%$  & $76.22\%$   & $78.04\%$  & $78.28\%$  & $78.32\%$  & \boldmath{$79.12\%$}& $62.86\%$  \\
\cline{2-9}
                        & $80\%$& $74.65\%$  & $74.90\%$   & $77.80\%$  & \boldmath{$78.11\%$}  & $76.74\%$  & $75.17\%$& $62.51\%$  \\
\cline{2-9}
                        & $85\%$& $71.01\%$  & $70.83\%$   & $75.69\%$  & \boldmath{$76.53\%$}  & $76.50\%$  & $73.00\%$& $60.73\%$  \\
\hline
\end{tabular}
\caption{Binary classification. Percentage of correctly classified samples in the testing set.}
\label{tab:bin_class_tab}
\end{table}


{\footnotesize
\begin{thebibliography}{10}

\bibitem{ding2019global}
X.~Ding, X.~Zhou, Y.~Guo, J.~Han, J.~Liu, {\em et~al.}, ``Global sparse
  momentum {SGD} for pruning very deep neural networks,'' {\em Advances in
  Neural Information Processing Systems}, vol.~32, 2019.

\bibitem{zimmer2022compression}
M.~Zimmer, C.~Spiegel, and S.~Pokutta, ``Compression-aware training of neural
  networks using {F}rank-{W}olfe,'' {\em arXiv preprint arXiv:2205.11921},
  2022.

\bibitem{duchi2011adaptive}
J.~Duchi, E.~Hazan, and Y.~Singer, ``Adaptive subgradient methods for online
  learning and stochastic optimization.,'' {\em Journal of machine learning
  research}, vol.~12, no.~7, 2011.

\bibitem{mcmahan2010adaptive}
H.~B. McMahan and M.~Streeter, ``Adaptive bound optimization for online convex
  optimization,'' {\em arXiv preprint arXiv:1002.4908}, 2010.

\bibitem{kingma2014adam}
D.~P. Kingma, ``{A}dam: a method for stochastic optimization,'' {\em arXiv
  preprint arXiv:1412.6980}, 2014.

\bibitem{tieleman2012lecture}
T.~Tieleman and G.~Hinton, ``Lecture 6.5-rmsprop, coursera: Neural networks for
  machine learning,'' {\em University of Toronto, Technical Report}, vol.~6,
  2012.

\bibitem{zeiler2012adadelta}
M.~D. Zeiler, ``{ADADELTA}: an adaptive learning rate method,'' {\em arXiv
  preprint arXiv:1212.5701}, 2012.

\bibitem{gratton2022first}
S.~Gratton, S.~Jerad, and P.~L. Toint, ``First-order objective-function-free
  optimization algorithms and their complexity,'' {\em arXiv preprint
  arXiv:2203.01757}, 2022.

\bibitem{gratton2023multilevel}
S.~Gratton, A.~Kopani{\v{c}}{\'a}kov{\'a}, and P.~L. Toint, ``Multilevel
  objective-function-free optimization with an application to neural networks
  training,'' {\em SIAM Journal on Optimization}, vol.~33, no.~4,
  pp.~2772--2800, 2023.

\bibitem{gratton2024complexity}
S.~Gratton, S.~Jerad, and P.~L. Toint, ``Complexity of a class of first-order
  objective-function-free optimization algorithms,'' {\em Optimization Methods
  and Software}, pp.~1--31, 2024.

\bibitem{conn2000trust}
A.~R. Conn, N.~I. Gould, and P.~L. Toint, {\em Trust region methods}.
\newblock SIAM, 2000.

\bibitem{yuan2015recent}
Y.-x. Yuan, ``Recent advances in trust region algorithms,'' {\em Mathematical
  Programming}, vol.~151, pp.~249--281, 2015.

\bibitem{reed1993pruning}
R.~Reed, ``Pruning algorithms-a survey,'' {\em IEEE transactions on Neural
  Networks}, vol.~4, no.~5, pp.~740--747, 1993.

\bibitem{lecun1989optimal}
Y.~LeCun, J.~Denker, and S.~Solla, ``Optimal brain damage,'' {\em Advances in
  neural information processing systems}, vol.~2, 1989.

\bibitem{zhu2017prune}
M.~Zhu and S.~Gupta, ``To prune, or not to prune: exploring the efficacy of
  pruning for model compression,'' {\em arXiv preprint arXiv:1710.01878}, 2017.

\bibitem{zhang2015accelerating}
X.~Zhang, J.~Zou, K.~He, and J.~Sun, ``Accelerating very deep convolutional
  networks for classification and detection,'' {\em IEEE transactions on
  pattern analysis and machine intelligence}, vol.~38, no.~10, pp.~1943--1955,
  2015.

\bibitem{denton2014exploiting}
E.~L. Denton, W.~Zaremba, J.~Bruna, Y.~LeCun, and R.~Fergus, ``Exploiting
  linear structure within convolutional networks for efficient evaluation,''
  {\em Advances in neural information processing systems}, vol.~27, 2014.

\bibitem{yu2017compressing}
X.~Yu, T.~Liu, X.~Wang, and D.~Tao, ``On compressing deep models by low rank
  and sparse decomposition,'' in {\em Proceedings of the IEEE conference on
  computer vision and pattern recognition}, pp.~7370--7379, 2017.

\bibitem{courbariaux2016binarized}
M.~Courbariaux, I.~Hubara, D.~Soudry, R.~El-Yaniv, and Y.~Bengio, ``Binarized
  neural networks: Training deep neural networks with weights and activations
  constrained to +1 or -1,'' {\em arXiv preprint arXiv:1602.02830}, 2016.

\bibitem{wang2018two}
P.~Wang, Q.~Hu, Y.~Zhang, C.~Zhang, Y.~Liu, and J.~Cheng, ``Two-step
  quantization for low-bit neural networks,'' in {\em Proceedings of the IEEE
  Conference on computer vision and pattern recognition}, pp.~4376--4384, 2018.

\bibitem{kim2020position}
J.~Kim, K.~Yoo, and N.~Kwak, ``Position-based scaled gradient for model
  quantization and pruning,'' {\em Advances in neural information processing
  systems}, vol.~33, pp.~20415--20426, 2020.

\bibitem{han2015learning}
S.~Han, J.~Pool, J.~Tran, and W.~Dally, ``Learning both weights and connections
  for efficient neural network,'' {\em Advances in neural information
  processing systems}, vol.~28, 2015.

\bibitem{yu2012exploiting}
D.~Yu, F.~Seide, G.~Li, and L.~Deng, ``Exploiting sparseness in deep neural
  networks for large vocabulary speech recognition,'' in {\em 2012 IEEE
  International conference on acoustics, speech and signal processing
  (ICASSP)}, pp.~4409--4412, IEEE, 2012.

\bibitem{louizos2017learning}
C.~Louizos, M.~Welling, and D.~P. Kingma, ``Learning sparse neural networks
  through $\ell_0$ regularization,'' {\em arXiv preprint arXiv:1712.01312},
  2017.

\bibitem{alvarez2017compression}
J.~M. Alvarez and M.~Salzmann, ``Compression-aware training of deep networks,''
  {\em Advances in neural information processing systems}, vol.~30, 2017.

\bibitem{hu2016network}
H.~Hu, ``Network trimming: A data-driven neuron pruning approach towards
  efficient deep architectures,'' {\em arXiv preprint arXiv:1607.03250}, 2016.

\bibitem{hoefler2021sparsity}
T.~Hoefler, D.~Alistarh, T.~Ben-Nun, N.~Dryden, and A.~Peste, ``Sparsity in
  deep learning: Pruning and growth for efficient inference and training in
  neural networks,'' {\em Journal of Machine Learning Research}, vol.~22,
  no.~241, pp.~1--124, 2021.

\bibitem{guo2016dynamic}
Y.~Guo, A.~Yao, and Y.~Chen, ``Dynamic network surgery for efficient {DNN}s,''
  {\em Advances in neural information processing systems}, vol.~29, 2016.

\bibitem{mocanu2018scalable}
D.~C. Mocanu, E.~Mocanu, P.~Stone, P.~H. Nguyen, M.~Gibescu, and A.~Liotta,
  ``Scalable training of artificial neural networks with adaptive sparse
  connectivity inspired by network science,'' {\em Nature communications},
  vol.~9, no.~1, p.~2383, 2018.

\bibitem{he2018soft}
Y.~He, G.~Kang, X.~Dong, Y.~Fu, and Y.~Yang, ``Soft filter pruning for
  accelerating deep convolutional neural networks,'' {\em arXiv preprint
  arXiv:1808.06866}, 2018.

\bibitem{pokutta2020deep}
S.~Pokutta, C.~Spiegel, and M.~Zimmer, ``Deep neural network training with
  {F}rank-{W}olfe,'' {\em arXiv preprint arXiv:2010.07243}, 2020.

\bibitem{molchanov2016pruning}
P.~Molchanov, S.~Tyree, T.~Karras, T.~Aila, and J.~Kautz, ``Pruning
  convolutional neural networks for resource efficient inference,'' {\em arXiv
  preprint arXiv:1611.06440}, 2016.

\bibitem{theis2018faster}
L.~Theis, I.~Korshunova, A.~Tejani, and F.~Husz{\'a}r, ``Faster gaze prediction
  with dense networks and {F}isher pruning,'' {\em arXiv preprint
  arXiv:1801.05787}, 2018.

\bibitem{lu2022learning}
M.~Lu, X.~Luo, T.~Chen, W.~Chen, D.~Liu, and Z.~Wang, ``Learning
  pruning-friendly networks via {F}rank-{W}olfe: One-shot, any-sparsity, and no
  retraining,'' in {\em International Conference on Learning Representations},
  2022.

\bibitem{argyriou2012sparse}
A.~Argyriou, R.~Foygel, and N.~Srebro, ``Sparse prediction with the $ k
  $-support norm,'' {\em Advances in Neural Information Processing Systems},
  vol.~25, 2012.

\bibitem{rao2017group}
N.~Rao, M.~Dud{\'\i}k, and Z.~Harchaoui, ``The group $k$-support norm for
  learning with structured sparsity,'' in {\em 2017 IEEE International
  Conference on Acoustics, Speech and Signal Processing (ICASSP)},
  pp.~2402--2406, IEEE, 2017.

\bibitem{frank1956algorithm}
M.~Frank, P.~Wolfe, {\em et~al.}, ``An algorithm for quadratic programming,''
  {\em Naval research logistics quarterly}, vol.~3, no.~1-2, pp.~95--110, 1956.

\bibitem{levitin1966constrained}
E.~S. Levitin and B.~T. Polyak, ``Constrained minimization methods,'' {\em USSR
  Computational mathematics and mathematical physics}, vol.~6, no.~5,
  pp.~1--50, 1966.

\bibitem{reddi2016stochastic}
S.~J. Reddi, S.~Sra, B.~P{\'o}czos, and A.~Smola, ``Stochastic {F}rank-{W}olfe
  methods for nonconvex optimization,'' in {\em 2016 54th annual Allerton
  conference on communication, control, and computing (Allerton)},
  pp.~1244--1251, IEEE, 2016.

\bibitem{WuWardBott18}
X.~Wu, R.~Ward, and L.~Bottou, ``Wngrad: Learn the learning rate in gradient
  descent,'' {\em arXiv preprint arXiv:1803.02865}, 2018.

\bibitem{Corletal96}
R.~M. Corless, G.~H. Gonnet, D.~E. Hare, D.~J. Jeffrey, and D.~E. Knuth, ``On
  the {L}ambert {W} function,'' {\em Advances in Computational mathematics},
  vol.~5, pp.~329--359, 1996.

\bibitem{van2007sparco}
E.~van~den Berg, M.~Friedlander, G.~Hennenfent, F.~Herrmann, R.~Saab, and
  O.~Y{\i}lmaz, ``Sparco: A testing framework for sparse reconstruction,'' {\em
  Dept. Comput. Sci., Univ. British Columbia, Vancouver, Tech. Rep.
  TR-2007-20,[Online]. Available: http://www. cs. ubc. ca/labs/scl/sparco},
  2007.

\bibitem{GratToin24}
S.~Gratton and P.~L. Toint, ``{S2MPJ} and {{\sf CUTEst}} optimization problems
  for {M}atlab, {P}ython and {J}ulia.'' arXiv:2407.07812, 2024.

\bibitem{wen2010fast}
Z.~Wen, W.~Yin, D.~Goldfarb, and Y.~Zhang, ``A fast algorithm for sparse
  reconstruction based on shrinkage, subspace optimization, and continuation,''
  {\em SIAM Journal on Scientific Computing}, vol.~32, no.~4, pp.~1832--1857,
  2010.

\bibitem{porcelli2014variable}
M.~Porcelli and F.~Rinaldi, ``A variable fixing version of the two-block
  nonlinear constrained {G}auss-{S}eidel algorithm for $\ell$ 1-regularized
  least-squares,'' {\em Computational Optimization and Applications}, vol.~59,
  no.~3, pp.~565--589, 2014.

\bibitem{ksvd}
M.~Aharon, M.~Elad, and A.~Bruckstein, ``{K-SVD}: An algorithm for designing
  overcomplete dictionaries for sparse representation,'' {\em IEEE Transactions
  on signal processing}, vol.~54, no.~11, pp.~4311--4322, 2006.

\bibitem{machrep}
$\_\_\_$, ``{UCI} machine learning repository,'' 2013.

\bibitem{chang2011libsvm}
C.-C. Chang and C.-J. Lin, ``{LIBSVM}: {A} library for support vector
  machines,'' {\em ACM transactions on intelligent systems and technology
  (TIST)}, vol.~2, no.~3, pp.~1--27, 2011.

\end{thebibliography}

}
\end{document}